# A relation between undrained CPTu results and the state parameter for liquefiable soils

Lluís Monforte [1*], Marcos Arroyo [1,2] and Antonio Gens [1,2]

1. Centre Internacional de Mètodes Numèrics en Enginyeria (CIMNE), Barcelona, Spain
2. Universitat Politècnica de Catalunya (UPC-BarcelonaTECH), Barcelona, Spain

* Centre Internacional de Mètodes Numèrics en Enginyeria (CIMNE)
Campus Nord UPC
C/ Gran Capità, s/n
08034 Barcelona (SPAIN)
lluis.monforte@upc.edu lmonforte@cimne.upc.edu

**Abstract:** This paper presents a new interpretation procedure to estimate the initial state parameter from cone penetration testing (CPTu) in undrained conditions based on the results from a comprehensive set of numerical simulations of CPTu in low permeability liquefiable soil. CPTu simulations are performed using the Particle Finite Element method, whereas the material response is modelled with CASM. The effects of soil constitutive parameters and roughness of the soil-steel interface are examined. It turns out that the numerical results are correctly summarized by an analytical relation derived from undrained cavity expansion results in critical state soils, as long as the constraints resulting from cone geometry are taken into account. The resulting adapted analytical formulation is notable for its simplicity and ease of use, comparing favorably with existing alternatives.

**Keywords:** CPTu test, Liquefaction, state parameter, tailings, PFEM

## 1. Introduction

Some geomaterials may exhibit a sudden loss of strength during undrained loading; this undrained softening is characterized by a decrease of the effective mean stress and a reduction of the mobilized undrained shear strength. The result of this sudden strength loss is known as static or flow liquefaction. The consequences of failures where static liquefaction is involved are momentous, as they give little or no warning of the impeding failure. A well-known example is the 2019 failure of tailings dam B1 in Brumadinho, Brazil (Arroyo & Gens, 2021).

Static liquefaction is best understood within the framework of critical state soil mechanics. Been and Jefferies (1985) introduced the concept of state parameter to characterize the behavior of soils. The state parameter is defined as the difference between the current void ratio and the void ratio on the Critical State Line at the same value of mean effective stress. A material will be dilatant if the state parameter is negative whereas if the state parameter is positive the material is contractive and may undergo undrained strength softening.

Knowing the state parameter of a soil will then be very helpful to establish if a soil is susceptible to flow liquefaction or not. Unfortunately, state-preserving extraction of samples from geomaterials susceptible to static liquefaction has proven very difficult (Been, 2016). If a measure of in situ state is necessary, it needs to be inferred from the results of in-situ tests. Amongst in situ tests, the cone penetration test with water pressure measurement (CPTu) enjoys a prominent position and much work has been done to infer state parameter values from CPTu records. These relations are generally based on normalized cone metrics. The notation employed in this work for normalized cone metrics is given in (Table 1). Following (Shuttle & Jefferies, 2016) normalization is based on the mean stress instead of the vertical stress.

The initial approach to derive state parameter values from CPTu measurements involved empirical correlation, using data acquired in calibration chambers (CC) with full-size CPTu (Jefferies & Been, 2016). A key result from this kind of work was that state parameter $\psi$ and normalized cone tip resistance $Q_p$ were related through the following expression (Been et al. 1987):

$$Q_p = \frac{(q_c - p_0)}{p_0'} = k \exp(-m\psi) \tag{1}$$

Where $q_c$ is cone tip resistance, $p_0$ and $p_0'$ represent, respectively, mean total and effective stress before testing and $k$ and $m$ were material dependent parameters. In particular, it was observed a strong influence of soil plastic compressibility as given by the critical state line slope in the compression plane, λ.

A disadvantage of this approach is that full-size CPTu in calibration chambers was practically limited to permeable materials, like sands. However, many soils susceptible to flow liquefaction are predominantly silt-sized and have relatively low permeabilities. Been et al (1989) addressed this problem. They observed that independently obtained empirical relations between $Q_p$ and overconsolidation ratio in clays could be cast in the same form of eq (1) simply by using a modified normalized tip strength that included pore pressure effects. This supported a new generalized form of the CPTu vs state parameter relation, given by

$$\frac{q_t - p_0}{p_0'}\left(1 - \frac{u_2 - u_0}{q_t - p_0}\right) = Q_p\left(1 - B_q\right) = \bar{k}\exp(-\bar{m}\psi) \tag{2}$$

where $B_q$ is the excess water pressure ratio and tip resistance, $q_t$ , is now corrected by unequal area effects. This relation was later developed by Plewes et al (1992) into a screening method for flow liquefaction, proposing empirical correlations for parameters $\bar{k}$ and $\bar{m}$. Plewes method carries significant uncertainty (Torres-Cruz, 2021) but offers reasonable central estimates (Reid, 2015) and is still frequently applied, particularly for tailings. Other empirical approaches to estimate state parameter are also in use (e.g. Robertson, 2009) but they are not applicable in undrained conditions.

In principle, modelling offers a way to obtain less uncertain relations between state parameter and cone metrics. Shuttle & Jefferies (1998) pioneered this approach using Nor-Sand, a critical state model that incorporates the state parameter concept. They solved numerically the problem of drained spherical cavity expansion. Their work showed that the limit pressures from cavity expansion had a relation with initial state parameter that, in a first approximation, could be also summarized by an exponential expression in the same form as Equation (1). Their numerical work also detected that the parameters $k$ and $m$ were not just a function of plastic compressibility λ, but also of other soil properties, notably normalized elastic stiffness, but also plastic hardening modulus and critical state friction.

Shuttle & Cunning (2007) applied the same procedure but now in a fully undrained condition. It was again observed, from spherical cavity expansion simulations with Nor-Sand, that the limiting cavity pressures obtained for a single set of constitutive parameters had a relatively simple relation with the state parameter. The expression used to summarize the results was very similar to that of the Plewes method, namely

$$\overline{Q_p}\left(1-\overline{B_q}\right)+1=\bar{k}\exp\left(-\bar{m}\psi\right) \tag{3}$$

where, following the convention of Table 1, the bars above the normalized metrics $Q_p$ and $B_q$ indicate that they correspond to cavity expansion. Shuttle & Jefferies (2016) revisited the problem, using an evolved version of Nor-Sand to simulate both drained and undrained spherical cavity expansion cases. They corroborated previous findings. Equation (3) was valid to summarize cavity expansion results for a specific drainage condition and a set of constitutive parameters. Parameters $\bar{k}$ and $\bar{m}$ were a complex function of various soil properties for the drained case, although depended only on critical state line properties (λ, M) for the undrained case.

Geometry and boundary conditions of expanding spherical cavities are very different from those in a CPTu. Even if the relation between limit cavity pressure and state parameter may be casted in the same form as that between normalized tip resistance and state parameter (e.g. Eq 3), the meaning of the variables employed is different (Table 1) and the values of the parameters $\bar{k}$ and $\bar{m}$ are also different. For instance, Shuttle & Jefferies (1998) obtained limit pressures from spherical cavity expansion in Ticino sand that, for the same state parameter, were one order of magnitude below the normalized cone tip resistance measured in CC tests.

These differences have been tackled by means of empirical shape correction functions, mapping cavity limit pressure into cone tip resistance. Shape correction functions were obtained comparing CC tests in sand with spherical cavity expansion results with NorSand by Ghafghazi & Shuttle (2008) and then revised by Shuttle & Jefferies (2016). Because CC testing was not available for undrained conditions the correction functions developed for tests in drained conditions are assumed valid for the undrained case.

Shape-corrected predictions of spherical cavity expansion simulations using site-specific calibrated NorSand model are currently considered state-of-the-art for state parameter evaluation. This has been facilitated by dissemination of accessible NorSand cavity expansion simulation tools ("widget"; Shuttle & Jefferies 2016). Nevertheless, direct validation for silts remains elusive due to complicating factors such as partial drainage (Reid & Smith, 2021) or sample disturbance (Fourie et al, 2022).

The cavity expansion method has always been conceived of (Shuttle & Jefferies, 1998, 2016) as an approximation, forced by the difficulties inherent to more realistic models of CPTu. Taking advantage of advances in numerical methods, CPTu testing can now be realistically simulated using different procedures (Lu et al, 2004; Nazem et al, 2012; Ceccato et al, 2017; Monforte et al, 2018a; Hauser and Schweiger, 2021). Numerical simulation provides new insights into the mechanisms that take place during the test, and can be used to assess the accuracy of current interpretation techniques (Monforte et al, 2018a, 2022a) or to propose enhancements (Monforte et al, 2018b). However, to tackle the problem of state parameter inversion in liquefiable materials it is also necessary to use a constitutive model that captures the essential mechanics at play.

Pezeshki & Ahmadi (2021) used a large deformation finite difference method to simulate CPTu in undrained conditions in a material represented by NorSand. Their results showed significant differences with those obtained using shape corrected cavity expansion and the same material description. They also compared the output of their parametric analyses with the empirical Plewes method and the differences were even more significant: for instance, they observed an effect of λ on the inversion parameter $\overline{k}$ opposite to that applied in Plewes method.

Monforte et al (2021) demonstrated how the Particle Finite Element Method (PFEM) could be employed to simulate CPTu in liquefiable materials represented with CASM (Clay and Sand Model). This constitutive model has been recently applied with success to analyse the onset of several flow liquefaction failures (Arroyo & Gens, 2021; Mánica et al. 2021, 2022). Monforte et al. (2022a) followed on that study to investigate the relation between cone results and undrained peak and residual strengths. These previous works did not focus on the relation of state parameter and cone metrics, which is addressed here.

The work presented in this paper is limited to cases in which cone penetration is undrained. The relative simplicity of this case allows to obtain closed form expressions of cavity expansion solutions that may be used as a reference to compare with the more precise results obtained with numerical simulation. This work is organized as follows: first, the adopted constitutive model is briefly summarized; then, a reference undrained cavity expansion solution is presented. After describing the numerical technique employed to carry out the simulations of CPTu testing, the results from a comprehensive CPTu simulation set are presented. Finally, a new interpretation technique is introduced.

## 2. Constitutive model

The constitutive model employed is a modified version of the isotropic, critical state-based model CASM, which was originally proposed by Yu (1998). CASM may represent a wide range of soil behaviors including softening, ductile or dilatant materials (Yu, 1998). A brief review of the constitutive model employed for the simulations of this study is presented here, more details can be found in Yu (1998) and Mánica el al. (2021; 2022).

In this work, CASM is applied within the framework of large-strain elasto-plasticity in which the deformation gradient splits multiplicatively into an elastic and a plastic part (Simo and Hughes, 1998). The elastic behavior is described by means of an hyperelastic law; the relation between the Kirchhoff stress tensor, $\boldsymbol{\tau}$, and Hencky elastic strains, $\boldsymbol{\epsilon}^e$, is expressed as (Houlsby, 1985):

$$\boldsymbol{\tau} = p_r \exp\left(-\frac{\epsilon_v^e}{\kappa^\star}\right) \mathbf{1} + 2\, G\, \boldsymbol{\epsilon}_d^e \tag{4}$$

where $\kappa^\star = \kappa / (1 + e_0)$, $\kappa$ is the slope of the reloading curve at the $\ln(p') - e$ space, $G$ is the shear modulus, $p_r$ is a reference pressure at which the void ratio is equal to $e_0$ and $\epsilon_v^e$ and $\boldsymbol{\epsilon}_d^e$ are the volumetric and deviatoric elastic Hencky strains. Thus, the bulk modulus depends linearly on the mean effective stress whereas the shear modulus is assumed constant, and it is computed from the initial mean effective stress as:

$$G = \frac{3\,(1 + 2\nu)}{2\,(1 + \nu)} \frac{p_0'}{\kappa^\star} \tag{5}$$

where $\nu$ is the Poisson's coefficient at the initial state and $p_0'$ is the initial mean effective stress.

The mathematical expression for the yield surface in triaxial compression is given by:

$$f = \left(\frac{q}{Mp'}\right)^n + \frac{1}{\ln r} \ln\left(\frac{p'}{p_c'}\right) \tag{6}$$

where $n$ and $r$ are model parameters controlling the shape of the yield surface, $p'_c$ is the preconsolidation stress, $p'$ is the Kirchhoff mean effective stress, $q$ is the deviatoric stress and $M$ is the slope of the CSL in $p' - q$ plane, related to the critical state friction angle, $\phi$, by

$$M = \frac{6 \sin \phi}{3 - \sin\phi} \tag{7}$$

For loading conditions different from triaxial compression the yield surface is generalized by having M depending on Lode-angle so as to obtain a smoothed Mohr-Coulomb envelope (Abbo and Sloan; 1995).

The classical isotropic volumetric hardening rule of critical state soil models is adopted:

$$p'_c = p_r \exp\left(-\frac{\epsilon_v^p}{\lambda^\star - \kappa^\star}\right) \tag{8}$$

where $\lambda^\star = \lambda / (1 + e_0)$, $\lambda$ is the slope of the Critical State line at the $\ln(p') - e$ space and $\epsilon_v^p$ is the Hencky volumetric plastic strain. The formulation of the constitutive model is closed by the flow rule proposed by Mánica et al (2021; 2022), whose dilatancy rule, for triaxial compression conditions, reads:

$$d^p = \frac{(m-1)}{m} \frac{(M^m - \eta^m)}{\eta^{m-1}} \tag{9}$$

where $m$ is a parameter of the model, $\eta = q/p'$ is the stress ratio, and $d^p$ is dilatancy, i.e. the ratio of incremental volumetric to shear plastic strain. The flow rule is generalized to non-triaxial conditions as done for the yield surface.

For the analytical solution presented in the next section, infinitesimal strains and small displacements are assumed; thus the elastic model and isotropic hardening law are formulated in terms of the Cauchy stress tensor and infinitesimal strains and the yield surface is expressed in terms of the Cauchy stress tensor. A flow rule is not specified, but it is assumed that during loading the soil tends to the critical state.

## 3. Reference cavity expansion solution

A closed form expression for the normalized effective resistance to undrained expansion of a cavity in CASM may be obtained if we assume that: (i) the soil at the cavity wall is at critical state and (ii) the initial state of the soil is isotropic and characterized by an initial mean effective stress $p'_0$ and initial state parameter $\psi_0$. To begin with, the effective limit cavity pressure, $\sigma'_c$, (i.e. the limit cavity pressure minus the water pressure) may be expressed in terms of the mean effective stress, $p'_{CS}$ , and deviatoric stress, $q_{CS}$, at critical state (Mo and Yu, 2017):

$$\sigma'_c = \sigma_c - u_c = \begin{cases} p'_{CS} + \frac{2q_{CS}}{3} = \left(1 + \frac{2M_\alpha}{3}\right) p'_{CS} & M_\alpha = \frac{6 \sin \phi}{3 - \sin \phi} \quad (Spherical) \\ p'_{CS} + \frac{q_{CS}}{2} = \left(1 + \frac{M_\alpha}{2}\right) p'_{CS} & M_\alpha = 2 \sin \phi \quad (Cylindrical) \end{cases} \tag{10}$$

Where $p'_{CS}$ and $q_{CS}$ stand for the mean effective stress and deviatoric stress at critical state after undrained loading from the initial stress state. $M_\alpha$ corresponds to the relevant CSL slope, which is that of triaxial compression for the spherical cavity case and of plane strain compression for the cylindrical case. Now, assuming undrained conditions, the mean effective stress and deviatoric stress at critical state are directly related to the initial mean effective stress, $p'_0$, the initial state parameter, $\psi_0$, and the slope of the critical state line in the $e - \ln p'$ plane, $\lambda$, as:

$$p'_{CS} = p'_0 \exp\left(\frac{-\psi_0}{\lambda}\right) \tag{11}$$

$$q_{cs} = M_\alpha \; p'_{cs} = p'_0 \, M_\alpha \; \exp\left(\frac{-\psi_0}{\lambda}\right)$$

Therefore, the effective cavity resistance may be expressed as:

$$\sigma'_c = (1 + \alpha M_\alpha)p'_0 \; \exp\left(\frac{-\psi_0}{\lambda}\right) \tag{12}$$

where $\alpha$ is 2/3 for the spherical case and 1/2 for the cylindrical one. This expression can be rewritten to show that the normalized effective cavity resistance (Table 1) is only dependent on the initial state parameter and two constitutive parameters describing the critical state line, $M_\alpha$ and $\lambda$:

$$\overline{Q_p}\left(1 - \overline{B_q}\right) + 1 = \frac{\sigma'_c}{p'_0} = (1 + \alpha M_\alpha) \exp\left(\frac{-\psi_0}{\lambda}\right) \tag{13}$$

Expressions for the other normalized cavity metrics (i.e. normalized resistance and excess water pressure ratio) may be also obtained using the analytical results of Mo and Yu (2017). However, the expressions obtained for these two metrics (see Appendix 1) depend not only on $\psi_0$, $\lambda$ and $M$ but also on the shear modulus $G$ and the yield surface shape parameters ($n$ and $r$).

Expression (13) has the same formal structure than the relation (Eq. 3) proposed by Shuttle & Cunning (2007). This is reasonable, as Eq. 3 was arrived at it through regression of numerical cavity expansion results on NorSand materials. In fact, the expression for normalized effective cavity resistance given above is not just specific of CASM, but it is also valid for any critical state model that uses the same description of the critical state line, such as NorSand. The main point, however, is that the cavity expansion solutions will offer an interesting comparison point for the more realistic numerical simulations that are presented next.

## 4. Numerical simulations

### *4.1 Overview*

We use the code G-PFEM to perform all CPTu simulations in this work. G-PFEM (Monforte et al 2017a, 2022b, Carbonell et al, 2022), has been developed for the analysis of large-strain problems in geomechanics, including the contact with rigid and deformable structures. G-PFEM –implemented in Kratos Multiphysics (Dadvand et al, 2010)- is based on PFEM, the Particle Finite Element Method (Idelsohn et al, 2004) whose key aspects include a Lagrangian description of motion, low order finite elements to compute the solution and the constant regeneration of the finite element mesh covering the domain.

Even if we are interested on behavior at the undrained limit, the problem is solved using effective stresses and a fully coupled hydromechanical formulation. The finite elements employed to discretize the domain are triangles with linear shape functions. These kinds of elements are known to be susceptible to volumetric locking; to alleviate this numerical pathology a mixed stabilized formulation is adopted, where nodal degrees of freedom include displacements and water pressure but also the local volume change (Jacobian) (Monforte et al, 2017b). The constitutive model is integrated using a classic scheme (Sloan et al, 2001) adapted to the large-strain formulation (Monforte et al, 2014). To increase robustness and reduce computational cost, the stress integration is formulated in the framework of the Implex technique (Oliver et al, 2008). To alleviate mesh-dependency of the solution induced by strain-softening, a nonlocal regularization technique is employed (Galavi & Schweiger, 2010).

The CPTu is simulated assuming axisymmetric conditions and a cone with standard dimensions. To avoid boundary interference, the cone is initially whished-in-place by 5.5D in a domain 20D wide and 30D high. The cone is then pushed into the soil at the standard velocity (0.02 m/s). The initial effective vertical stress is 100 kPa and the horizontal stress is given by $K_0 = 0.6$ (calculated with Jacky's formula using the CS friction angle). The initial vertical stress is maintained at the top of the domain, null vertical displacements are prescribed at the bottom and null radial displacements are set at the axis of symmetry and the outer radial boundary. A hydraulic conductivity of $k = 10^{-10}$ m/s is adopted in all simulations, which makes the analyses effectively undrained (see Monforte et al. 2021, for a parametric study on permeability).

An initial simulation series was performed using a smooth CPTu and five different materials (Table 2). These five materials are all slightly overconsolidated (isotropic OCR = 1.1) They share an isotropic compression line with $\lambda$ = 0.054 passing through $p_r$ = 100 kPa for $e_0$ = 1. They also share most constitutive parameters including elastic properties, ($\kappa = 0.016$, $\nu \approx 0.33$), critical state friction, ($\phi =$ 25 ), and flow-rule (m = 2.5). However, their Critical State Lines in the $\ln p' - e$ plane have different intercepts, $\Gamma_c$ . As a result, their initial state parameter is also different, but always positive, as we are interested in contractive materials. The constitutive parameters (and those employed in the parametric analysis) are not representative of any soil in particular, but are characteristic of a broad range of materials that might undergo static liquefaction (Torres-Cruz, 2019; Tarragó, 2021; Macedo and Vergaray, 2021).

To characterize the five reference materials, Figure 1 reports the simulation of undrained triaxial compression from their initial state. The peak undrained shear strength is almost coincident for all cases, whereas they all have different residual strengths. At the end of the simulation, all tests end at critical state, both at the $p' - q$ and $\ln p' - e$ planes. The values of undrained peak and residual strength as well as that of Bishop brittleness index are reported in Table 2.

After the initial series, parametric studies explored the effect of different factors in the results. Those studies always included 5 materials with the same characteristics as in the initial series except for the parameter that is varied. The results presented below include parametric studies to the effect of different soil properties and interface friction at the cone.

### *4.2 Initial simulation series*

Figure 2 shows the results of the CPTu simulations for the initial series. The basic results for this series include $q_n$, net tip resistance (tip resistance minus vertical total stress) and the excess pore pressures measured at mid-face ($\Delta u_1$) and cone shoulder ($\Delta u_2$) positions. The normalized metrics presented include normalized cone tip resistance, $Q_p$, and the normalized excess pore pressure at the cone shoulder, $B_{q2}$, and at the cone face, $B_{q1}$. A stationary state for all these metrics is observed after a penetration of 10 radii. Both net tip resistance and excess pore water pressure increase as the initial state parameter decreases, but the change induced in pore water pressure is much smaller. As a result, the increase in normalized tip resistance with decreasing state parameter results in a decrease of the normalized excess water pressure ratio.

Steady state normalized metrics for the numerical simulations are collected in Table 3 and presented in terms of the initial state parameter in Figure 3; for completeness, results employing water pressure measured at both reading positions are included. Also included in this figure are the values obtained with the reference infinite cavity expansion solutions, employing the same constitutive parameters. For a given material, total cavity resistance and excess pore pressure are evaluated using the formulae reported in Appendix 1 whereas normalized effective resistance is computed with Eq. (13). The resulting datapoints are then joined with a spline.

The effect of initial state parameter on normalized resistance and excess pore pressure is qualitatively similar for the CPTu G-PFEM simulation and for the cavity expansion solutions but, in general, there are significant quantitative differences between results from the cone simulations and those from cavity expansion. As the material parameters are the same, the observed differences are due to the geometrical and boundary condition differences between the two problems. Cone tip resistance is severely underestimated by cavity expansion solutions (Shuttle & Jefferies, 1998; Salgado & Prezzi, 2007), although they are closer to cone measurements when excess pore pressure is considered (Silva et al. 2006). The presence of tip resistance in both $Q_p$ and $B_q$ metrics largely explains the observed differences in these metrics between the G-PFEM and the cavity expansion solutions. This is not the case for the normalized effective tip resistance, where cavity expansion results show a remarkable quantitative agreement with the cone simulation when pore pressure measurements at the position $u_1$ are used to evaluate it.

To understand this result is necessary to examine in detail other simulation outputs. Figure 4 presents paths of invariant stress and excess pore pressure observed, as the cone approaches, at a point initially located 17 radii below the tip of the cone and 1.15 radii off the axis of symmetry. Results are presented for the most brittle material, with the largest initial state parameter (Material A) and for the less brittle material, with the smallest state parameter (Material E). The undrained cone advance pushes the material around the cone towards critical state. When the state parameter is initially high attaining critical state involves a large decrease in mean effective stress and, consequently, on shear stress. This is not the case when the initial state parameter is smaller, because critical state is attained without much stress reduction. In both cases the excess pore pressure is very similar, with large pore pressures just above the tip of the cone.

Figure 5 shows in more detail the stress state of the soil around the cone tip once steady state conditions have been achieved (further penetration does not modify the cone resistance nor the stress state around the cone tip). The near-liquefied condition that the cone induces in Material A is visible in the almost nil mean effective stress around the cone. It is also visible how the pore pressure field would have some radial symmetry if observed from the $u_1$ position but much less so if observed from the $u_2$ position, as the pore pressures along the shaft are very different from those along the tip face. This explains the greater similitude to cavity expansion results of cone metrics based on $u_1$

Figure 4 also shows that, for both materials, the soil remains close to triaxial compression conditions (Lode angle $\theta_L \approx$ -30º) and the effective stress path coincides with that of an undrained compression triaxial test, with some oscillations due to the non-local stress integration technique and the constant remeshing of PFEM. For more detail on stress paths during cone penetration see Monforte et al (2021).

To further characterize the stress state around the cone tip and to compare the effective stress field with that of a spherical cavity expansion, Figure 6 reports the magnitude and direction of the first principal effective stress for the most (A) and least (E) brittle materials. In the first column the first principal stress induced by CPTu penetration is normalized by the radial effective stress at critical state predicted by the expansion of a spherical cavity, Equation (10). Around the tip, the cone induces the same principal stresses than spherical cavity expansion, but the area over which this happens gets much closer to the tip in the more brittle material. This result explains the differences in normalized tip resistance observed in Figure 3a, where the absolute difference with cavity expansion increases with $\psi_0$.

For the initial simulation series, where the cone is considered smooth, the major principal effective stress is normal to the cone tip. This is clearly visible in the stress field plots (middle column) and contours of principal direction angle with the vertical shown in Figure 6. The initial state parameter of the different materials has no effect on this result.

### 4.3 *Effect of soil properties*

Plastic compressibility was the first factor identified by previous research (e.g. Been et al. 1987) as influential in the relation of cone metrics and state parameter. Two additional sets of simulations examine this effect: in the first, $\lambda$ is increased -from 0.056 to 0.106- while maintaining the same $\kappa$ value as in the initial series; in the second set $\lambda$ is also increased, but the plastic volumetric ratio ($\Lambda = (\lambda-\kappa) / \lambda$) is fixed. Because the initial void ratio and stress state is maintained, the increased values of $\lambda$ result in larger values of the state parameter for all the examined materials (Table 3).

Figure 7 presents the results of this parametric study. There is a major effect of $\lambda$ on the normalized cone metrics. On the other hand cone metrics are almost indifferent to the value of $\kappa$. As plastic compressibility increases the influence of initial state on cone metrics reduces. The cavity expansion solutions also show this effect in their normalized metrics, but the difference with CPTu simulations is still large, except, again, for normalized effective tip resistance when computed with pore pressure measured at the $u_1$ position.

Figure 8 shows the results of a parametric study to critical state friction. The effect of increasing the friction angle of the soil from 25º ($M$ = 0.98) to 33º ($M$ = 1.33), while changing $K_0$ from 0.6 to 0.47, as per Jaky's formula, is to shift up significantly the normalized cone resistance and slightly reduce the excess water pressure ratio. The combined effect of those changes on normalized effective resistance appears moderate, particularly if the metric is computed with pore pressure measured at the $u_1$position: when increases on $Q_p(1 - B_{q1}) + 1$ are around 10%. The changes due to critical state friction in the cavity expansion reference solutions are consistent with those observed in the simulation. Note that for the normalized effective resistance the (small) effect of friction (going from M = 1.0 to 1.33) can be readily derived from eq. (13).

The significant role that elastic stiffness might play on the relation between state parameter and CPTu metrics was only appreciated at a relatively late stage, when cavity expansion solutions for NorSand became available (Shuttle & Jefferies, 1998). The role of elastic stiffness is here examined by a simulation series in which Poisson coefficient changed from 0.33 to 0.2. All the other things being equal, this implies an almost doubling of the elastic shear modulus, $G$, as per Eq (5). Such doubling will also affect other metrics usually employed in this context – like the rigidity index, $I_r = G/S_u$, or the dimensionless rigidity $G/p'$, (Jefferies et al 2016).

The results of the parametric study on elastic stiffness are also plotted in Figure 8. In this particular case, the change in elastic stiffness $G$ had almost the same effect on $Q_p$ than the change on critical state friction. Elastic stiffness contributes to increased tip resistance by increasing containment of the plasticized zone around the cone. Referring to Figure 6 (left column) the high stress region away from the tip is, essentially, an elastic zone. Elastic parameters do not feature in the reference cavity expansion solution for normalized effective resistance, Equation (13); it is noticeable that their effect on normalized effective resistance (Figure 8 c) appears to be also very small in the simulation results.

### *4.4 Effect of friction at the cone interface*

In all the previous analyses, a smooth cone-soil interface behavior has been assumed. However, friction mobilized at the soil-cone interface in CPTu tests will depend on cone roughness, which is generally poorly controlled. Currently, CPTu friction sleeves are manufactured with average roughness, $R_a$, between 0.65μm and 0.15μm (EN ISO 22476-1). For fine grained soils the resulting normalized roughness ($R_a$ / $D_{50}$) will typically vary in the range ($10^{-3}$ to $10^{-1}$). This range of normalized roughness leads to interface efficiencies between 0.3 and 0.9 (Subba Rao et al. 2000). On this basis, the higher values of interface efficiency would be representative for the finer, more impermeable, range of liquefiable soils for which undrained cone penetration is relevant.

To explore this effect, we have repeated the simulations assuming a rough interface, with a contact friction angle of 18º, implying a mobilized interface friction efficiency of 0.73. The results are presented in Figure 9; they show that higher friction at the interface increases a little the normalized tip resistance and decreases significantly the normalized excess water pressures. These effects are larger for materials with lower state parameters. The effective tip resistance also increases with friction; employing the $u_1$ measurement, simulation results plot slightly above the predictions for a spherical cavity whereas those using the $u_2$ reading plot higher, with a slope similar to that predicted by cavity expansion.

Figure 6 shows that the stress state around the rough cone has some significant differences respect to the smooth case. These differences are not so much in the magnitude of the principal stress (left column) but rather on its orientation (middle and right columns). For instance, the direction of major principal stress at the tip of the cone has an angle around 150º with respect to the vertical, and is no longer normal to the cone tip. Also, there is a consistent orientation of principal stress along the shaft instead of the more rapidly changing and disordered picture that emerges for the smooth case.

## 5. A simplified inversion formula for undrained conditions

### *5.1 Development of the inversion formula*

A remarkably consistent result from the simulations is that the normalized effective cone resistance, if computed with pore pressure measured at the $u_1$ position, is much closer to the analogue cavity expansion results than it is the case for other metrics. This suggests that a useful inversion formula to obtain state parameter from the cone may be obtained based on the cavity expansion expression (12). There is only a small difference between normalized effective cone resistance computed for the spherical and cylindrical cavity solutions -the difference between $\alpha = 0.5$ and $\alpha = 0.66$ in Equation (12). However, the stress fields around the tip, particularly that of pore pressure around the $u_1$ position, seem closer to spherical symmetry (Figure 5). Besides, joint examination of all the comparisons presented above between cavity expansion and numerical results suggest that the spherical case performs slightly better. For these reasons we will use the spherical cavity solution to develop the simplified inversion formula.

To do so consider first the ratio between normalized effective cavity resistance and the corresponding cone metric, which is

$$\frac{Q_p\left(1-B_{q1}\right)+1}{\overline{Q_p}\left(1-\overline{B_q}\right)+1}=\frac{q_c-u_1}{\sigma_c'}=c_q \tag{14}$$

The ratio $c_q$ can be conceived as a geometric correction factor, similar to those applied in previous work (Shuttle & Cunning, 2007). An analytical expression for $c_q$ may be obtained under some simplifying assumptions about the stress field near the tip of the cone. The first one is that the stress field is homogeneous and in conditions corresponding to triaxial compression. The second one is that the soil is at critical state. Lastly, we assume that the pore pressure measured at the $u_1$ position is representative of the homogeneous pore pressure acting on the face. The numerical results previously discussed (Figure 4 to Figure 6) indicate that these hypotheses are approximately correct, with the last one likely to imply a larger error. As shown in Appendix 2, these hypotheses lead to:

$$c_q = \frac{q_c-u_1}{\sigma_c'} = \frac{M+3M\cos(2\rho)-3\sqrt{3}\,M\sin(2\rho)+6}{4M+6} \tag{15}$$

where $\rho$ is the angle between the vertical and the major principal direction of the stress tensor. This expression is evaluated in Figure 10 for different values of $M$. It has a value of 1 for smooth cones (when $\rho = 120^{\circ}$) and larger values for rough cones ($\rho > 120^{\circ}$). Therefore, the value of $c_q$ not only depends on soil properties but also on those of the soil-cone interface.

The cone excess water pressure ratio in Equation (14) is $B_{q1}$, evaluated with excess water pressure measured at the $u_1$ position. In practice, pore pressure measurement at the $u_2$ position are more generally available than those at $u_1$ (Lunne et al, 2022). Results from the G-PFEM simulations show (Figure 11) that the ratio $\beta = B_{q1}/B_{q2}$ is not very sensitive to the state parameter or to soil properties, although there is, again, some effect of cone roughness, with higher friction at the interface resulting in slightly increased $\beta$. It is worth noting that the values of $\beta$ observed in the simulations are well aligned with field observations in normally consolidated soils (Sills et al, 1989; Peuchen et al, 2010).

Summarizing, the proposed interpretation technique to estimate the state parameter from cone readings is expressed as:

$$Q_p\left(1-B_{q1}\right)+1 = Q_p\left(1-\beta\,B_{q2}\right)+1 = c_q\left(1+\frac{2}{3}M\right)\exp\left(-\frac{\psi_0}{\lambda}\right) \tag{16}$$

*5.2 Comparison with previous expressions*

The inversion formula (16) has a similar structure to the classical inversion formula (2), if we take as inversion parameters $\overline{k} = c_q(1 + \alpha M)$ and $\overline{m} = 1/\lambda$. It is interesting to examine the effect of λ on these inversion parameters according to the formula proposed here and compare it with previous proposals by Plewes et al. (1992) and Pezeshki and Ahmadi (2021). This comparison is only indicative, as each proposal uses a slightly different definition of the normalized effective tip resistance (Table 4 ).

The comparison is presented in Figure 12 for $M = 1.4$, a fairly typical value for tailings (Macedo & Vergaray, 2022). For the range of λ that is more likely to be encountered in practice – say $\lambda = 0.05$ to 0.2 - the different proposals for $\overline{m}$ are close at the high-compressibility end, but they diverge for low-compressibility materials. The higher values of $\overline{m}$ that the new method predicts in that situation will result in smaller state parameters, other things being equal.

There are also significant differences for $\overline{k}$ , a parameter that corresponds to the normalized effective tip resistance for a soil that is already at critical state (i.e. has $\psi_0 = 0$). In Plewes et al (1992) method $\overline{k}$ continuously decreases as λ increases. Pezeshki and Ahmadi (2021) predict the opposite effect - in fact the proposal even predicts negative normalized effective tip stresses for $\lambda < 0.0106$. Our proposal implies that $\overline{k}$ is independent of λ and will mostly depend on critical state friction (M).

The different inversion equations discussed may be checked using the simulation results (Figure 13). The simplified method proposed here performs well, which is not surprising, since the simplified method was inspired by the simulation results. The Plewes method always overestimates the input state parameter of the simulations. The performance of the Pezeshki-Ahmadi method is reasonable for low values of initial state parameter, but significantly overestimates the input value when $\psi_0$ increases.

It will be highly desirable to check this and other methods with well-controlled test results. Field CPTu records might be exploited for this purpose, but in most cases the uncertainty derived from poor sample quality is a significant obstacle to direct validation. The best hope for this appears to lie in newer scaled-down calibration chambers -Sadrekarimi & Jones, 2022; Ayala et al., 2022- that have been recently perfected to test finer-grained materials. Currently, however, that database remains limited, and, for instance, it does not include results for CPTu advanced in undrained or partly drained conditions.

## 6. Discussion

*6.1 Applicability of the undrained case*

Liquefaction is a behaviour characteristic of low plasticity, predominantly granular soils, like sands and silts. It may be then questioned how relevant is an analysis based on undrained CPTu response. In fact, even if the simulations above have been run using a very low value of hydraulic conductivity ($k_h$ = $10^{-10}$ m/s), an undrained response to the CPTu appears in far more permeable materials.

Using G-PFEM simulation of CPTu in Cam-Clay soils, Monforte et al (2018b) showed that hydraulic conductivity values below $10^{-6}$ m/s do not allow enough flow of water around the cone during penetration to significantly modify the measured pore pressure. This condition may be used to characterize undrained CPTu advance. Monforte et al (2018b) also showed how this result was coherent with well-accepted backbone curves (De Jong & Randolph, 2012) describing the transition from drained to undrained CPTu behavior. Several important cases of known flow liquefaction involve materials with values of $k_h$ below this limit. For instance, the fine tailings in Brumadinho were characterized with $k_h$ < $10^{-7}$ m/s (Robertson et al. 2019; Arroyo & Gens, 2021). The tailings close to the failed section of Merriespruit dam had measured $k_h$ < $2\text{x}10^{-7}$ m/s (Stryrdom & Williams, 1999). The hydraulic fill that liquefied Prat quay in Barcelona harbor (Gens, 2022) had an estimated hydraulic conductivity of $5\text{x}10^{-8}$ m/s (Tarrago, 2021)

Discussing liquefaction in tailings from an applied perspective, Fourie et al. (2022) modified the chart of Schneider et al. (2008) to suggest that when $B_q$ > 0.2 the CPTu test might be considered as undrained. More than 50% of the fine tailings in Brumadinho had $B_q$ > 0.2, and $B_q$ was systematically above 1 in the vicinity of the borehole section that is likely to have triggered the failure (Arroyo & Gens, 2021). Very high values of $B_q$ were also registered in the liquefiable materials of Cadia (Jefferies et al. 2019). Many other instances of potentially liquefiable silty tailings of high $B_q$ may be found in the literature (e.g. Shuttle & Jefferies, 2016; dos Santos Junior et al. 2022; Bonin et al. 2022). It seems then that even a solution restricted to undrained CPTu advance is likely to find application in the characterization of liquefiable silt*s*.

*6.2 Effect of $K_0$ on the estimated state parameter*

The use of mean effective stress to normalize cone metrics has the inconvenient of requiring an independent estimate of the coefficient of earth-pressure at rest $K_0$. This is a parameter that is difficult to estimate; this difficulty will introduce uncertainty in the estimations of state parameter.

To see how important this effect might be, note that the equation relating CPTu measurements and state parameter may be inverted to give

$$\psi = \left(-\frac{1}{\overline{m}}\right)(\ln Q' - \ln \overline{k}) = -\frac{\ln Q'}{\overline{m}} + \frac{\ln \overline{k}}{\overline{m}} \tag{17}$$

where:

$$Q' = Q_p\left(1 - \beta B_{q2}\right) + 1 = \frac{q_c - \beta \Delta u_2 - u_0}{p_0'} = \frac{q_c - \beta \Delta u_2 - u_0}{\sigma_{v0}'} \frac{3}{1 + 2K_0} \tag{18}$$

leading to

$$\psi = -\frac{1}{\overline{m}} \ln\left(\frac{q_c - \beta \Delta u_2 - u_0}{\sigma_{v0}' \overline{k}}\right) + \frac{1}{\overline{m}} \ln\left(\frac{3}{1 + 2K_0}\right) = \psi_{iso} + \Delta\psi_{K_0} \tag{19}$$

where the inferred state parameter is expressed as the sum of $\psi_{iso}$, the value that would be obtained in isotropic stress conditions (i.e. when $K_0 = 1$) and $\Delta\psi_{K_0}$, a correction term due to the actual value of $K_0$. The correction term only depends on $K_0$ and, through $\overline{m}$, on the plastic compressibility λ.

Liquefiable silts are typically normally consolidated or slightly overconsolidated and thus $K_0$ is likely to vary within a relatively narrow range, between 0.5 and 1. For instance, Arroyo and Gens (2021) employed a value of 0.5 in the numerical analysis of Dam B1 of Brumadinho and $K_0 = 0.6$ in the analysis of Merriespruit (Mánica et al, 2022), whereas values of 0.7 were used in the interpretation of the CPTu records of Rose Creek (Shuttle & Cunning, 2007) and Neves-Corvo (Shuttle & Jefferies, 2016). A value of $K_0$ = 0.7 was also used in the numerical simulation of Cadia Dam failure (Jefferies et al. 2019).

Figure 14 illustrates the value of the correction term due to $K_0$ as a function of the plastic compressibility λ for the method presented in this work. If no information is available about the value of $K_0$ application of the method may proceed assuming that $K_0$ = 0.7. It appears that the error in state parameter estimation that may be incurred in that case will be small (less than 0,02), except for the more compressible materials.

## 7. Conclusions

This work has reported a comprehensive parametric analysis of CPTu testing in undrained, liquefiable materials. The effect of the initial state parameter, most of the constitutive parameters of the soil (volumetric compressibility, shear modulus, and parameters describing the critical state line) and the roughness of the soil-cone interface have been examined. The G-PFEM simulations have clarified several of the mechanisms affecting the relation between cone metrics and initial state parameter in undrained conditions, as well as the limits of cavity expansion as an analogue model of CPTu advance.

G-PFEM may be applied in further studies on state parameter inversion. Perhaps the first question to address should be on the effect of partial or full drainage on the relation between cone metrics and state parameter. Other questions of relevance may involve more elaborate constitutive descriptions, for instance the effect of anisotropy (in both hydraulic and mechanical behavior) or of the curvature of the CSL. The computational efficiency of G-PFEM offers an attractive platform for this kind of study.

Based on the G-PFEM results a simplified analytical interpretation technique to estimate the initial state parameter from CPTu results has been proposed. The new proposal clarifies the origin of correction factors applied to relate cavity expansion and cone results. The results shown emphasize the interpretation advantages that would follow from measuring pore pressure at the cone face. They also suggest that more attention should be paid to measuring cone interface friction to improve the reliability of inversion procedures. There are significant differences between the predictions of this new method and those of previous methods for state inversion from undrained CPTu. More work is needed to investigate the origin of these differences and the reliability of the different methods.

## 8. APPENDIX 1: Additional relation for cavity expansion in CASM

Mo and Yu (2017) report cavity expansion solution for CASM in undrained conditions. The main hypotheses are that the soil is homogeneous, and the initial stress state is isotropic, with null water pressure, initial mean effective stress of $p_0'$ and initial preconsolidation of $p_{c0}$. In the case of an infinite expansion, the limit cavity stress can be expressed as (Mo & Yu, 2017):

$$\sigma_c = p_0' - \frac{m_d}{m_d + 1} q_{cs} \ln A_3 + 2\, G_0 m\, A_4 \tag{20}$$

where $m_d$ is equal to 1 for cylindrical cavity and 2 for a spherical cavity, $q_{cs}$ is the deviatoric stress at critical state (Equation (11)), $G_0$ stands for the shear modulus and:

$$A_3 = 1 - \exp\left[-\left(\frac{\ln(R_0)}{\ln(r)}\right)^{\frac{1}{n}} \frac{M_\alpha\, p_0'}{2\, G_0}\right] \tag{21}$$

$$A_4 = \frac{1}{m_d + 1} \sum_{k=1}^{\infty} \frac{A_3^k}{k^2} \tag{22}$$

where $R_0 = p_{c0}/p_0'$ is the ratio between the initial preconsolidation stress to mean effective stress.

The effective cavity resistance can be expressed as (see Equation (6)):

$$\sigma_c' = \left(1 + \frac{m_d M_\alpha}{m_d + 1}\right) p_{cs}' \tag{23}$$

being $p_{cs}'$ the mean effective stress at critical state after undrained triaxial loading from the initial state.

The water pressure at the cavity can be computed as the difference between the total and the effective cavity resistance:

$$u_c = \sigma_c - \sigma_c' \tag{24}$$

## 9. APPENDIX 2: Derivation of the geometrical term

This appendix presents a complete derivation of Equation (15). We assume that the effective stress state close to the cone tip is homogeneous, at critical state and in triaxial compression conditions. If the major principal stress has an angle $\rho$ with respect to the vertical (Figure 15), the total stress tensor is given by:

$$\begin{aligned} \boldsymbol{\sigma} &= \begin{bmatrix} \sigma_r & \sigma_{rz} & 0 \\ \sigma_{rz} & \sigma_z & 0 \\ 0 & 0 & \sigma_\theta \end{bmatrix} \\ &= \begin{bmatrix} \cos(\rho) & \sin(\rho) & 0 \\ -\sin(\rho) & \cos(\rho) & 0 \\ 0 & 0 & 1 \end{bmatrix} \cdot \begin{bmatrix} \sigma_3' & 0 & 0 \\ 0 & \sigma_1' & 0 \\ 0 & 0 & \sigma_3' \end{bmatrix} \cdot \begin{bmatrix} \cos(\rho) & -\sin(\rho) & 0 \\ \sin(\rho) & \cos(\rho) & 0 \\ 0 & 0 & 1 \end{bmatrix} \\ &+ \begin{bmatrix} u & 0 & 0 \\ 0 & u & 0 \\ 0 & 0 & u \end{bmatrix} \end{aligned} \tag{25}$$

Where the $\sigma_i'$ represent the principal stresses and $u$ stands for the water pressure. In triaxial compression at the critical state $q_{cs} = M\, p_{cs}'$, and the principal stresses are $\sigma_1' = p_{cs}' + \frac{2q_{cs}}{3} = p_{cs}'\left(1 + M\frac{2}{3}\right)$ and $\sigma_2' = \sigma_3' = p_{cs}' - \frac{q_{cs}}{3} = p_{cs}'\left(1 - \frac{M}{3}\right)$.

The traction vector acting at the face of the cone is:

$$\boldsymbol{t} = \boldsymbol{\sigma} \cdot \boldsymbol{n} \tag{26}$$

where $\boldsymbol{n}$ is the normal to the cone. As the stress state is assumed homogenous, the cone total resistance may be obtained as:

$$q_c = t_z \frac{\pi R^2 / \cos(60)}{\pi R^2} \tag{27}$$

where $t_z$ is the vertical component of $\boldsymbol{t}$ and the last term is the ratio between the cone face area and the projected area of the cone.

Using Matlab symbolic computations, the following expression may be obtained for the cone tip resistance:

$$q_c = p'_{cs} + p'_{cs} M \left( \frac{1}{6} + \frac{\cos(2\rho)}{2} - \frac{\sqrt{3} \sin(2\rho)}{2} \right) + u \tag{28}$$

Finally, by subtracting the water pressure and dividing by $\sigma'_c$ the cavity stress at critical state in spherical condition – given by Equation (10)- we obtain:

$$\frac{q_c - u}{\sigma'_c} = \frac{M + 3 M \cos(2\rho) - 3\sqrt{3} \sin(2\rho) + 6}{4M + 6} \tag{29}$$

If we suppose now that the water pressure that would be measured mid-face (i.e. at the $u_1$ position) is representative of the water pressure at the face we recover Equation (15).

## 10. Acknowledgements

The authors wish to thank Dr. Stefano Collico (UPC), Mr. Davide Besenzon (Escuela Superior Politécnica del Litoral) and Mohammad Razavinasab (CIMNE) for fruitful discussions at the beginning of the research.

Financial support of Ministerio de Ciencia e Innovación of Spain (MCIN/AEI/10.13039/501100011033) through the Severo Ochoa Centre of Excellence project (CEX2018-000797-S) and research project PID2020-119598RB-I00 is gratefully appreciated.

## 11. Competing interests

The authors declared that there is no conflict of interest.

## 12. Data availability

Data generated or analyzed during this study are available from the corresponding author upon reasonable request.

## 13. NOTATION

| | |
|---|---|
| $A_3, A_4$ | Auxiliary variables. Appendix 1 |
| $B_q$ | Excess water pressure ratio |
| $B_{q1}$ | Excess water pressure ratio using $u_1$ reading |
| $B_{q2}$ | Excess water pressure ratio using $u_2$ reading |
| $\overline{B_q}$ | Cavity excess water pressure ratio |
| $c_q$ | Ratio between normalized effective cavity resistance and normalized effective cone resistance |
| $d^p$ | Dilatancy |
| $D_{50}$ | Median particle size |
| $e$ | Void ratio |
| $e_0$ | Initial void ratio |
| $G$ | Shear modulus |
| $G_0$ | Reference shear modulus. Appendix 1 |
| $K_0$ | Coefficient of earth pressure at rest |
| $k$ | Permeability |
| $k$ | Inversion parameter. Drained conditions. Equation (1) |
| $\overline{k}$ | Inversion parameter. Undrained conditions. Equation (2) |
| $M$ | Slope of the CSL in $p' - q$ plane in triaxial compression conditions |
| $M_\alpha$ | Slope of the CSL in $p' - q$ plane for different loading conditions. Equation (10) |
| $m$ | Constitutive parameter controlling the dilatancy. Equation (9) |
| $m$ | Inversion parameter. Drained conditions. Equation (1) |
| $m_d$ | Auxiliary variable. Appendix 1 |
| $\overline{m}$ | Inversion parameter. Undrained conditions. Equation (2) |
| | |
| $n$ | Shape parameter of CASM yield surface. Equation (6) |
| $\boldsymbol{n}$ | Normal |
| $p_r$ | Reference pressure. Equations (4) and (8) |
| $p'_c$ | Preconsolidation stress |
| $p'$ | Mean effective stress |
| $p'_0$ | Initial mean effective stress |
| $p'_{cs}$ | Mean effective stress at critical state conditions |
| $p_0$ | Initial mean total stress |
| $q$ | Deviatoric stress |
| $q_{cs}$ | Deviatoric stress at critical state conditions |
| $q_c$ | Cone tip resistance |
| $q_t$ | Cone tip resistance corrected by unequal area effects. |
| $Q_p$ | Normalized cone tip resistance |
| $\overline{Q_p}$ | Normalized cavity resistance |
| $r$ | Spacing ratio of CASM yield surface. Equation (6) |
| $r$ | Radial coordinate |
| $R$ | Cone radius |
| $R_a$ | Average roughness |
| $R_0 = p_{c0}/p'_0$ | Ratio between the initial preconsolidation stress to mean effective stress |
| $S_u$ | Undrained shear strenght |
| $\boldsymbol{t}$ | Traction vector |
| $u_0$ | Initial pore pressure |
| $u_1$ | Pore pressure at the $u_1$ position |
| $u_2$ | Pore pressure at the $u_2$ position |
| $u_c$ | Pore pressure at the cavity in limit condition |
| $z$ | Vertical coordinate |
| | |
| $\alpha$ | Parameter for the cavity expansion solution. Equation (13) |
| $\beta$ | Ratio between $\Delta u_2$ to $\Delta u_1$ |
| $\Gamma_c$ | Position of the CSL |
| $\delta$ | Soil-cone interface friction angle |
| $\boldsymbol{\epsilon}^{\mathrm{e}}$ | Elastic Hencky strain tensor |
| $\epsilon_v^{\mathrm{e}}$ | Volumetric elastic Hencky strain |
| $\boldsymbol{\epsilon}_d^{\mathrm{e}}$ | Deviatoric elastic Hencky strain |
| $\epsilon_v^{\mathrm{p}}$ | Volumetric plastic Hencky strain |
| $\eta = \frac{p'}{q}$ | Stress ratio |
| $\theta$ | Circumferential coordinate |
| $\theta_L$ | Lode's Angle |
| $\kappa^\star = \frac{k}{1+e_0}$ | Elastic compressibility |
| $\kappa$ | Slope of the reloading curve at the $\ln(p') - e$ |
| $\lambda^\star = \frac{\lambda}{1+e_0}$ | Plastic compressibility |
| $\lambda$ | Plastic compressibility at the $\ln(p') - e$ |
| $\Lambda = \frac{\lambda - \kappa}{\lambda}$ | Plastic volumetric ratio |
| $\nu$ | Poisson's coefficient |
| $\rho$ | Angle of the major principal stress with respect to the vertical |
| $\boldsymbol{\sigma}$ | Total Cauchy stress tensor |
| $\sigma'_1$ | Major principal effective stress |
| $\sigma'_3$ | Minor principal effective stress |
| $\sigma_c$ | Limit cavity resistance |
| $\sigma'_c$ | Limit effective cavity resistance |
| $\boldsymbol{\tau}$ | Kirchhoff stress tensor |
| $\phi$ | Critical state friction angle |
| $\psi$ | State parameter |
| $\psi_0$ | Initial state parameter |
| $\psi_{iso}$ | Inverted state parameter assuming isotropic conditions ($K_0 = 1$) |
| $\Delta\psi_{K_0}$ | Correction of the inverted state parameters due to $K_0$ |

## 15. TABLES

*Table 1: Definition of normalized cavity and cone metrics.*

| | Cavity expansion | Cone |
|---|---|---|
| Basic variables | $\sigma_c$ limit cavity pressure<br>$u_c$ pore pressure at the cavity in limit condition | $q_c$ cone tip resistance<br>$u_x$ pore pressure at transducer location "x" |
| Normalized resistance | $\overline{Q_p} = \frac{\sigma_c - p_0}{p_0'}$ | $Q_p = \frac{q_c - p_0}{p_0'}$ |
| Excess water pressure ratio | $\overline{B_q} = \frac{u_c - u_0}{\sigma_c - p_0}$ | $B_{qx} = \frac{u_x - u_0}{q_c - p_0}$ |
| Normalized effective resistance (Houlsby, 1989) | $\overline{Q_p}(1 - \overline{B_q}) + 1 = \frac{\sigma_c'}{p_0'} = \frac{\sigma_c - u_c}{p_0'}$ | $Q_p(1 - B_{qx}) + 1 = \frac{q_c - u_x}{p_0'}$ |

*Table 2: Constitutive parameters describing the shape of the yield function of the reference series. In addition,* $\kappa = 0.016$, $\lambda = 0.054$, *m = 2.5,* $e_0 = 1$ *at* 100 *kPa,* $\phi = 25$ , $\nu \approx 0.33$ *and OCR = 1.1*

| Material | n | r | $\Gamma_c$ (100 kPa) | $\psi_0$ | $S_u^p$ (kPa) | $S_u^{res}$ (kPa) | $I_b$ |
|---|---|---|---|---|---|---|---|
| A | 10 | 12 | 0.908 | 0.089 | 26.12 | 6.78 | 0.741 |
| B | 10 | 8 | 0.923 | 0.074 | 26.55 | 8.98 | 0.662 |
| C | 9 | 6 | 0.934 | 0.063 | 26.31 | 11.02 | 0.581 |
| D | 7 | 3 | 0.960 | 0.037 | 26.47 | 18.03 | 0.319 |
| E | 4 | 2 | 0.974 | 0.020 | 26.70 | 24.93 | 0.066 |

*Table 3: Material characteristics and CPTu and cavity expansion results from the different simulation series.*

| Series | Material | $\psi_0$ | $S_u^{peak}$ (kPa) | $S_u^{res}$ (kPa) | $Q_p$ | $B_{q1}$ | $B_{q2}$ | $Q_p(1-B_{q1})+1$ | $Q_p(1-B_{q2})+1$ | $Q_p(1-1.2B_{q2})+1$ | $\overline{Q_p}(1-\overline{B_q})+1$ Spherical | $\overline{Q_p}(1-\overline{B_q})+1$ Cylindrical |
|---|---|---|---|---|---|---|---|---|---|---|---|---|
| Reference | A | 0.0887 | 26.12 | 6.78 | 1.93 | 1.37 | 1.15 | 0.29 | 0.71 | 0.27 | 0.31 | 0.27 |
| | B | 0.0738 | 26.55 | 8.98 | 2.17 | 1.32 | 1.06 | 0.32 | 0.86 | 0.40 | 0.41 | 0.35 |
| | C | 0.0629 | 26.31 | 11.02 | 2.31 | 1.25 | 1.05 | 0.42 | 0.89 | 0.40 | 0.51 | 0.43 |
| | D | 0.0368 | 26.47 | 18.03 | 2.94 | 1.10 | 0.91 | 0.69 | 1.28 | 0.74 | 0.83 | 0.71 |
| | E | 0.0196 | 26.70 | 24.93 | 3.72 | 0.98 | 0.82 | 1.06 | 1.68 | 1.07 | 1.14 | 0.98 |
| Rough | A | 0.0887 | 26.12 | 6.78 | 2.10 | 1.30 | 0.95 | 0.36 | 1.11 | 0.71 | 0.31 | 0.27 |
| $\delta = 19º$ | B | 0.0738 | 26.55 | 8.98 | 2.36 | 1.25 | 0.99 | 0.40 | 1.01 | 0.54 | 0.41 | 0.35 |
| | C | 0.0629 | 26.31 | 11.02 | 2.55 | 1.19 | 0.90 | 0.50 | 1.25 | 0.79 | 0.51 | 0.43 |
| | D | 0.0368 | 26.47 | 18.03 | 3.45 | 1.01 | 0.70 | 0.98 | 2.04 | 1.56 | 0.83 | 0.71 |
| | E | 0.0196 | 26.70 | 24.93 | 4.20 | 0.92 | 0.59 | 1.34 | 2.71 | 2.22 | 1.14 | 0.98 |
| $\lambda$ | A | 0.1774 | 26.12 | 6.78 | 1.57 | 1.42 | 1.12 | 0.33 | 0.82 | 0.47 | 0.31 | 0.27 |
| $\lambda$ = 0.106 | B | 0.1475 | 26.55 | 8.98 | 1.77 | 1.33 | 1.07 | 0.41 | 0.88 | 0.50 | 0.41 | 0.35 |
| $\kappa$ = 0.032 | C | 0.1259 | 26.31 | 11.02 | 2.07 | 1.23 | 1.02 | 0.52 | 0.96 | 0.54 | 0.51 | 0.43 |
| | D | 0.0736 | 26.47 | 18.03 | 2.53 | 1.09 | 0.88 | 0.78 | 1.30 | 0.85 | 0.83 | 0.71 |
| | E | 0.0393 | 26.70 | 24.93 | 3.22 | 0.98 | 0.81 | 1.06 | 1.62 | 1.10 | 1.14 | 0.98 |
| Plastic ratio. $\Lambda$ | A | 0.2157 | 25.96 | 4.72 | 1.43 | 1.52 | 1.20 | 0.25 | 0.72 | 0.37 | 0.22 | 0.19 |
| $\lambda$ = 0.106 | B | 0.1794 | 26.37 | 6.64 | 1.71 | 1.43 | 1.15 | 0.27 | 0.75 | 0.36 | 0.30 | 0.26 |
| $\kappa$ = 0.016) | C | 0.1531 | 26.10 | 8.52 | 1.97 | 1.34 | 1.11 | 0.34 | 0.79 | 0.35 | 0.39 | 0.34 |
| | D | 0.0896 | 26.14 | 15.51 | 2.44 | 1.13 | 0.90 | 0.67 | 1.25 | 0.81 | 0.71 | 0.61 |
| | E | 0.0478 | 26.06 | 23.01 | 3.36 | 1.00 | 0.80 | 0.99 | 1.68 | 1.15 | 1.05 | 0.91 |
| $M$ | A | 0.0882 | 31.70 | 8.22 | 2.20 | 1.33 | 1.12 | 0.27 | 0.73 | 0.23 | 0.36 | 0.29 |
| $\phi$= 33º | B | 0.0734 | 32.16 | 10.88 | 2.55 | 1.24 | 1.03 | 0.40 | 0.93 | 0.41 | 0.47 | 0.39 |
| | C | 0.0625 | 31.95 | 13.37 | 2.77 | 1.19 | 1.00 | 0.48 | 1.01 | 0.46 | 0.58 | 0.48 |
| | D | 0.0362 | 32.22 | 21.95 | 3.68 | 1.04 | 0.86 | 0.85 | 1.51 | 0.87 | 0.95 | 0.78 |
| | E | 0.0187 | 32.73 | 30.55 | 4.65 | 0.98 | 0.81 | 1.07 | 1.87 | 1.12 | 1.33 | 1.09 |
| $\nu$ | A | 0.0887 | 26.12 | 6.78 | 2.30 | 1.33 | 1.13 | 0.25 | 0.69 | 0.17 | 0.31 | 0.27 |
| $\nu = 0.2$ | B | 0.0738 | 26.54 | 8.98 | 2.54 | 1.26 | 1.06 | 0.33 | 0.86 | 0.32 | 0.41 | 0.35 |
| | C | 0.0629 | 26.31 | 11.02 | 2.72 | 1.20 | 0.98 | 0.45 | 1.06 | 0.53 | 0.51 | 0.43 |
| | D | 0.0368 | 26.47 | 18.03 | 3.45 | 1.06 | 0.85 | 0.80 | 1.51 | 0.92 | 0.83 | 0.71 |
| | E | 0.0196 | 26.70 | 24.93 | 4.31 | 0.91 | 0.73 | 1.38 | 2.18 | 1.55 | 1.14 | 0.98 |

*Table 4: Interpretation techniques of the state parameter from CPTu testing in undrained conditions*

| | Plewes | Pezeshki and Ahmadi | This work |
|---|---|---|---|
| Normalized resistance | $Q_p = \frac{q_c - p_0}{p_0'}$ | $Q_p = \frac{q_c - p_0}{p_0'}$ | $Q_p = \frac{q_c - p_0}{p_0'}$ |
| Excess water pressure ratio | $B_q = \frac{\Delta u_2}{q_c - p_0}$ | $B_q = \frac{\Delta u_2}{q_c - \sigma_{v0}}$ | $B_q = \frac{\Delta u_1}{q_c - p_0} \approx \frac{\beta\, \Delta u_2}{q_c - p_0}$ |
| $\overline{k}$ | $M\left(3 + \frac{0.37}{\lambda}\right)$ | $M\left(3.3 - \frac{0.035}{\lambda}\right)$ | $c_q\left(1 + \frac{2}{3}M\right)$ |
| $\overline{m}$ | $11.9 - 30.62\lambda$ | $6 + \frac{0.1735}{\lambda}$ | $\frac{1}{\lambda}$ |

## 16. FIGURES

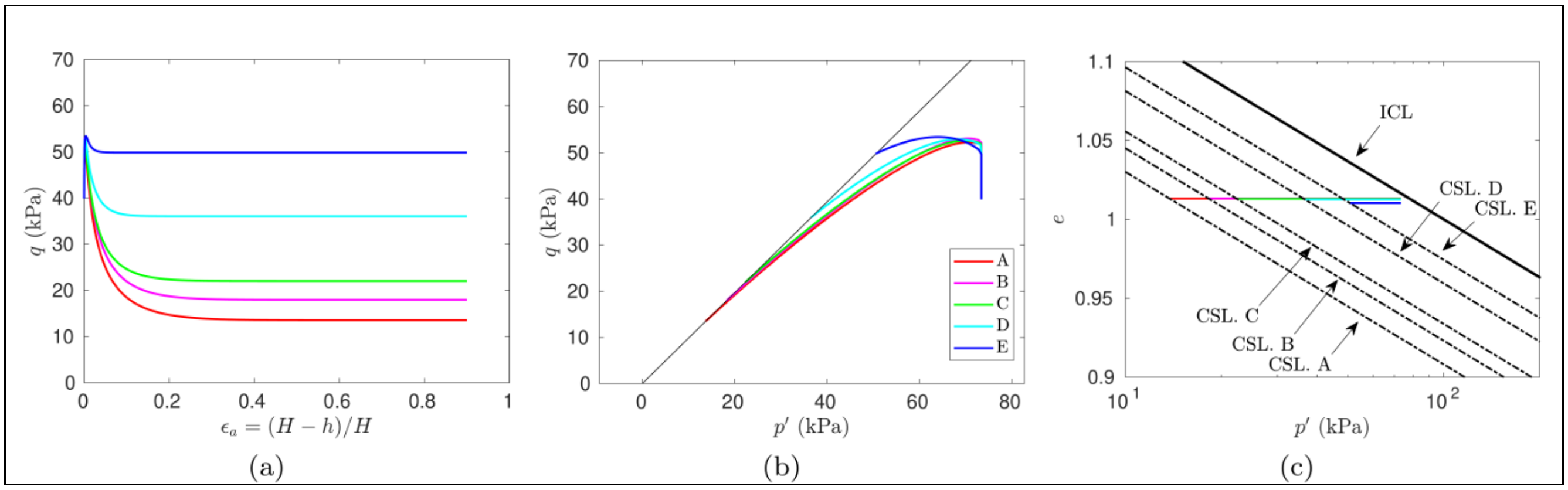

*Figure 1: Undrained triaxial response for the 5 materials in the initial simulation series (a) axial strain vs deviatoric stress (b) stress paths (c) trajectory in the compression plane, with indication of the different critical state lines of the materials.*

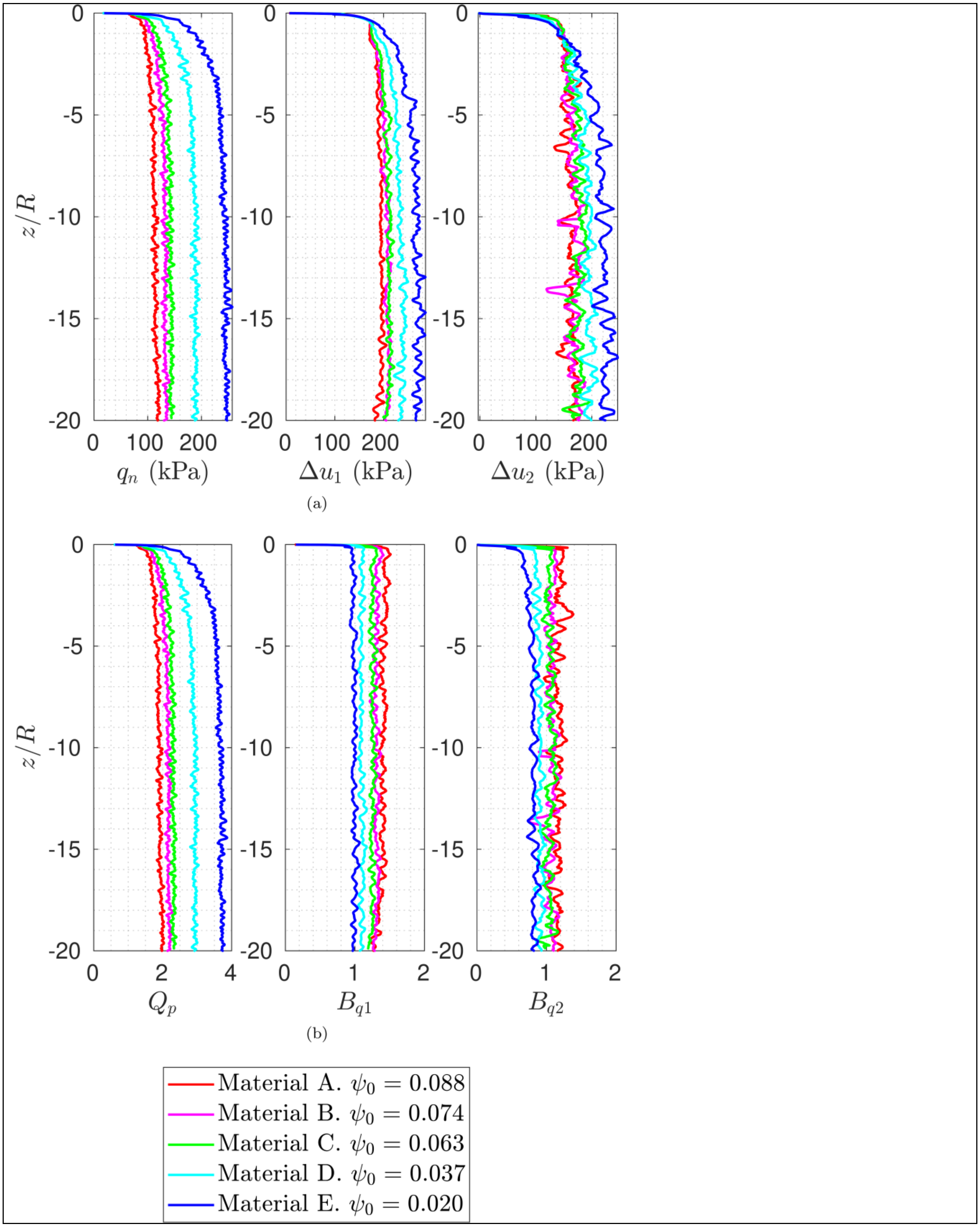

Figure 2: Initial CPTu simulation series. Basic measured variables (a) and normalized metrics (b).

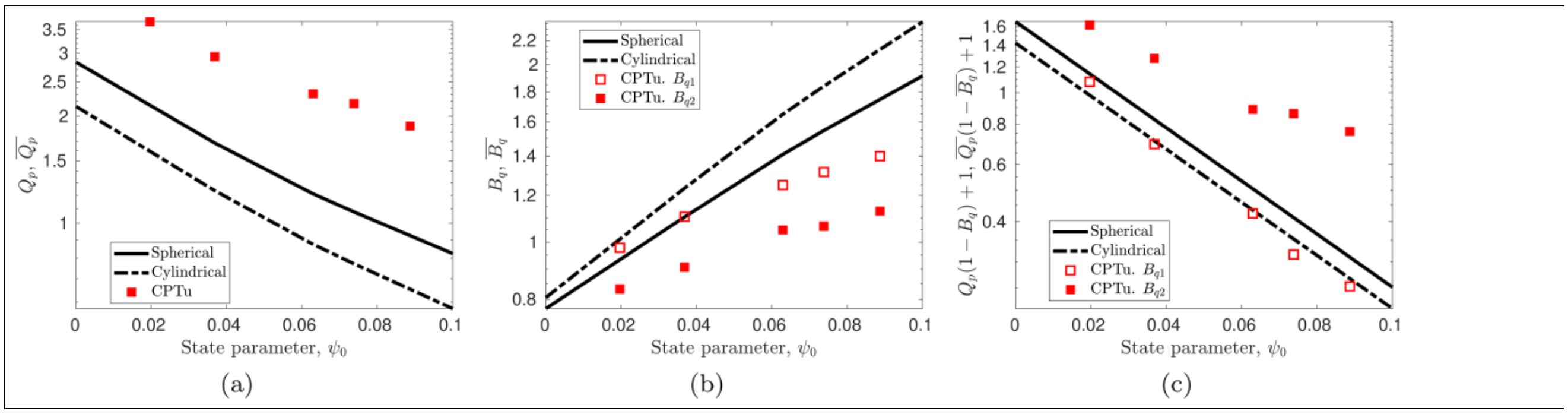

*Figure 3: Initial simulation series. Comparison of CPTu simulation results and cavity expansion reference solutions. Effect of initial state parameter on normalized resistance (a) excess water pressure ratio (b) and normalized effective resistance (c).*

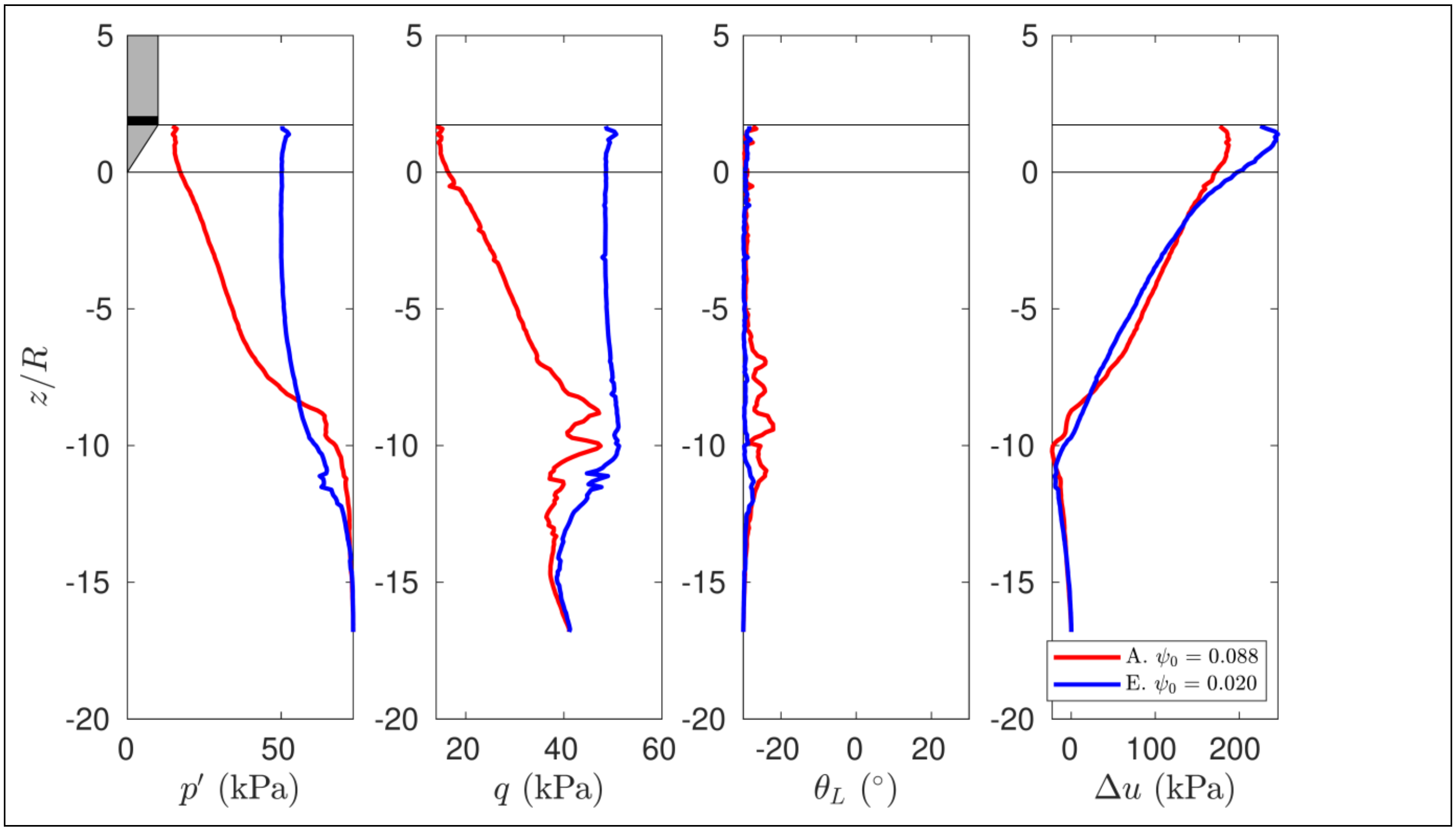

*Figure 4: Initial simulations. Stress paths (mean effective stress, deviatoric stress, Lode's angle and excess pore pressure) observed, as the cone tip approaches, at a point initially 17R below the cone tip and 1.15R off axis. The vertical axis represents normalized vertical offset between cone tip and observation point; values are negative when the cone tip is above the observation point. The cases represented correspond to the largest (A) and smaller initial state parameter (E).*

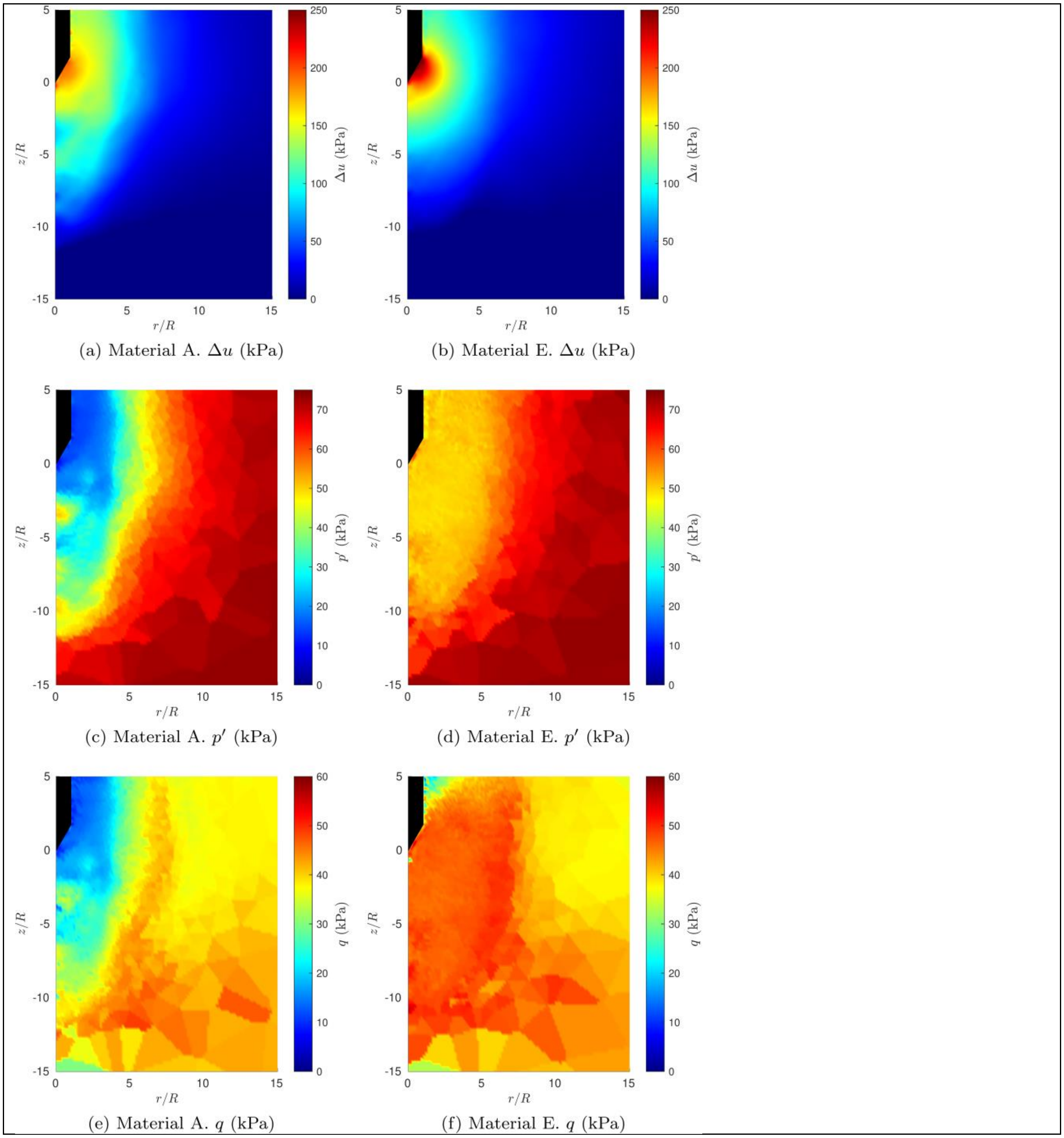

*Figure 5: Initial simulation series. Steady state contour plots of water pressure, (a) and (b), mean effective stress, (c) and (d), and deviatoric stresses, (e) and (f), for the case with largest (Material A, left) and smallest ( Material E, right) initial state parameter.*

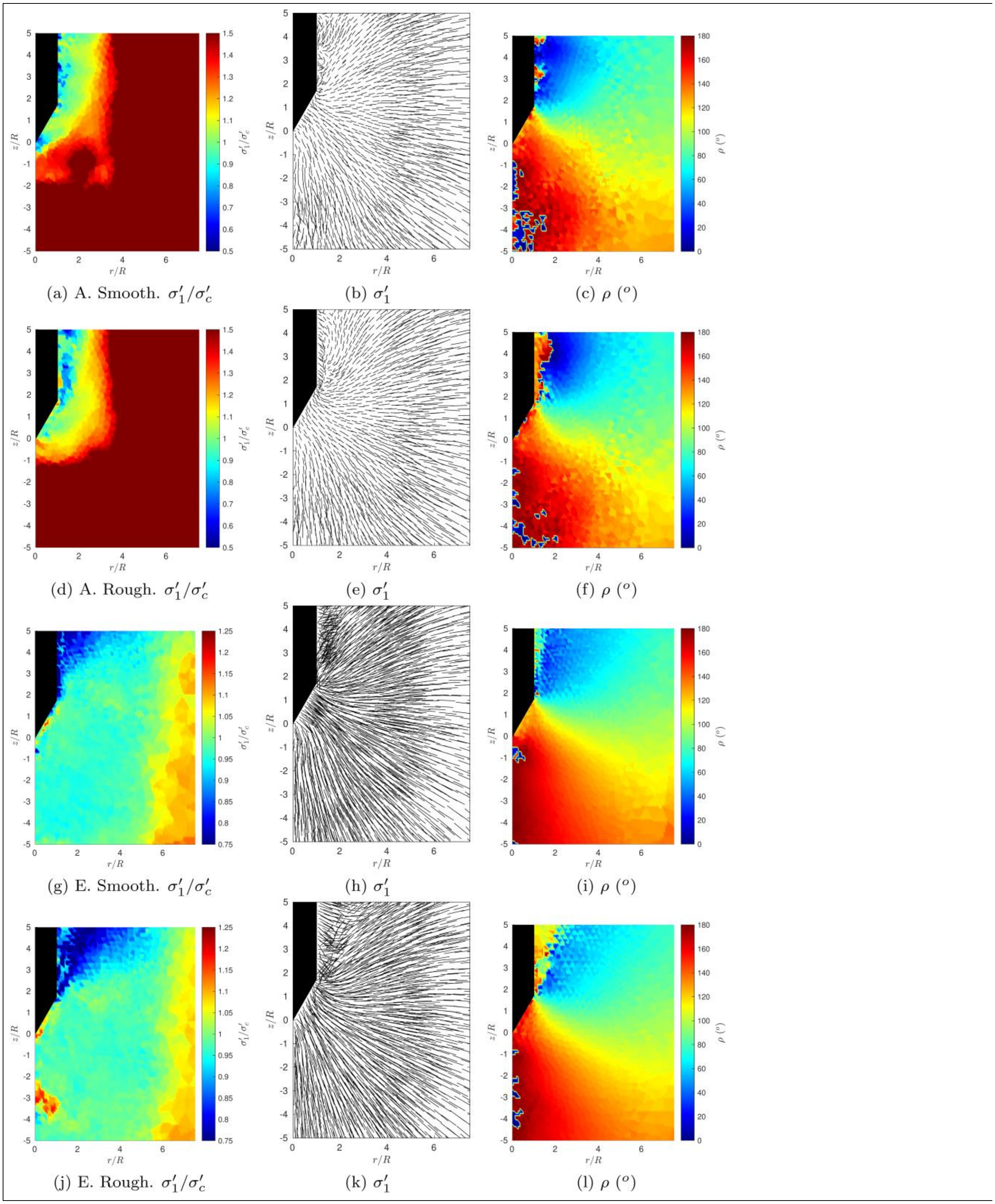

*Figure 6:Maximum principal effective stress normalized by the effective cavity resistance (left), maximum principal effective stress (center), direction respect the vertical of the maximum principal effective normalized stress (right). Material A, smooth (top), Material A, rough (second row), Material E, smooth (third row), Material E, rough (bottom).*

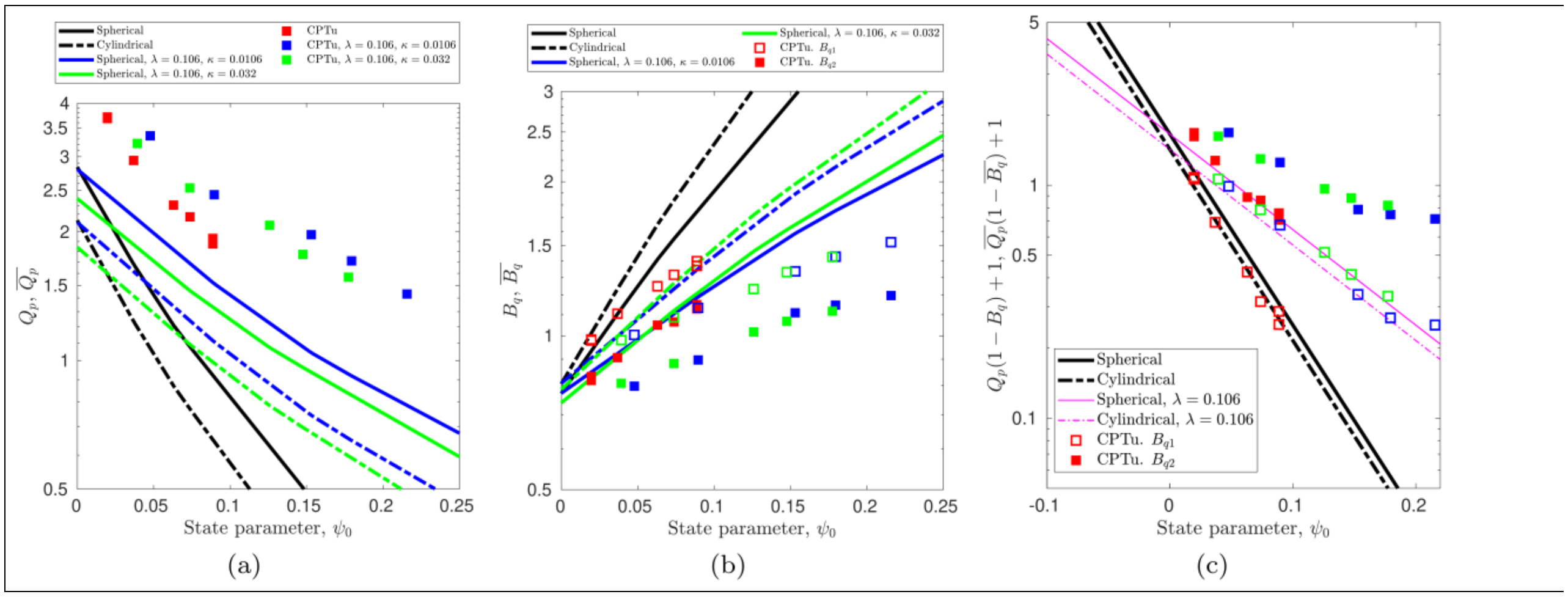

*Figure 7: CPTu simulation results vs cavity expansion reference solutions. Effect of initial state parameter on normalized resistance (a) excess water pressure ratio (b) and normalized effective resistance (c). Parametric study on soil plastic compressibility.*

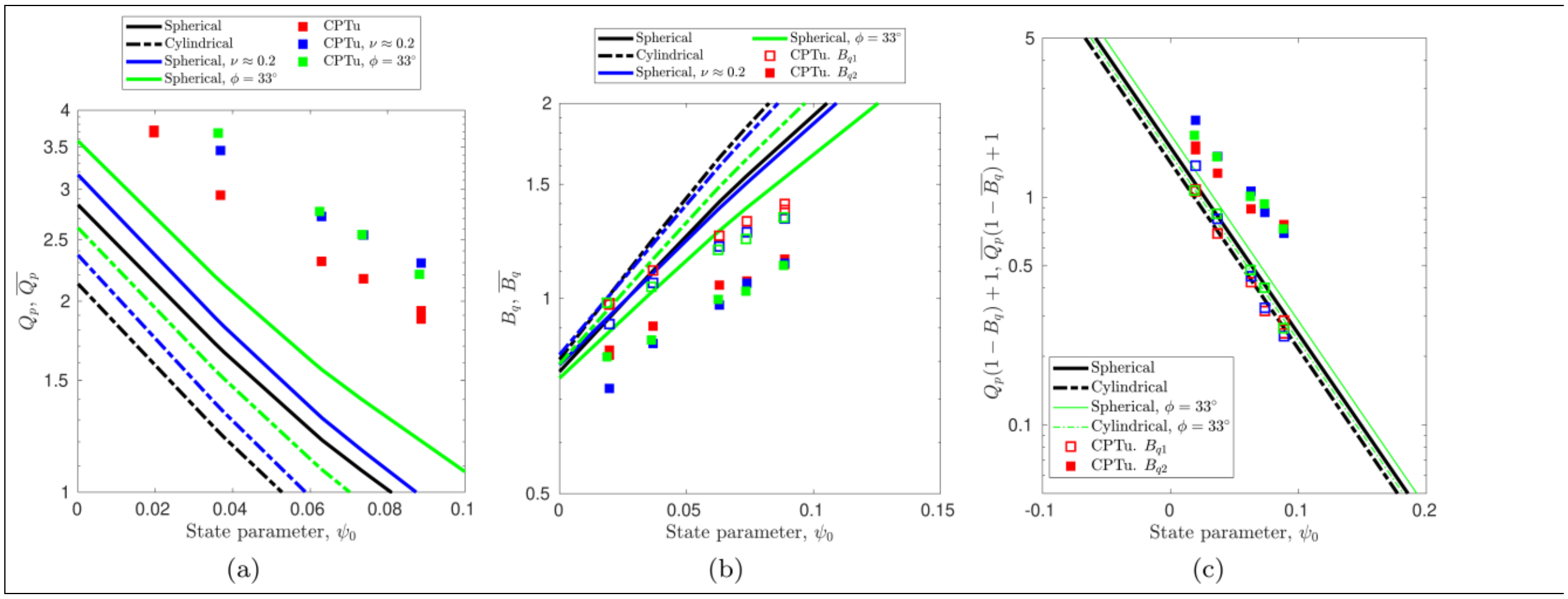

*Figure 8: CPTu simulation results vs cavity expansion reference solutions. Effect of initial state parameter on normalized resistance (a) excess water pressure ratio (b) and normalized effective resistance (c). Parametric study on soil friction and elastic stiffness.*

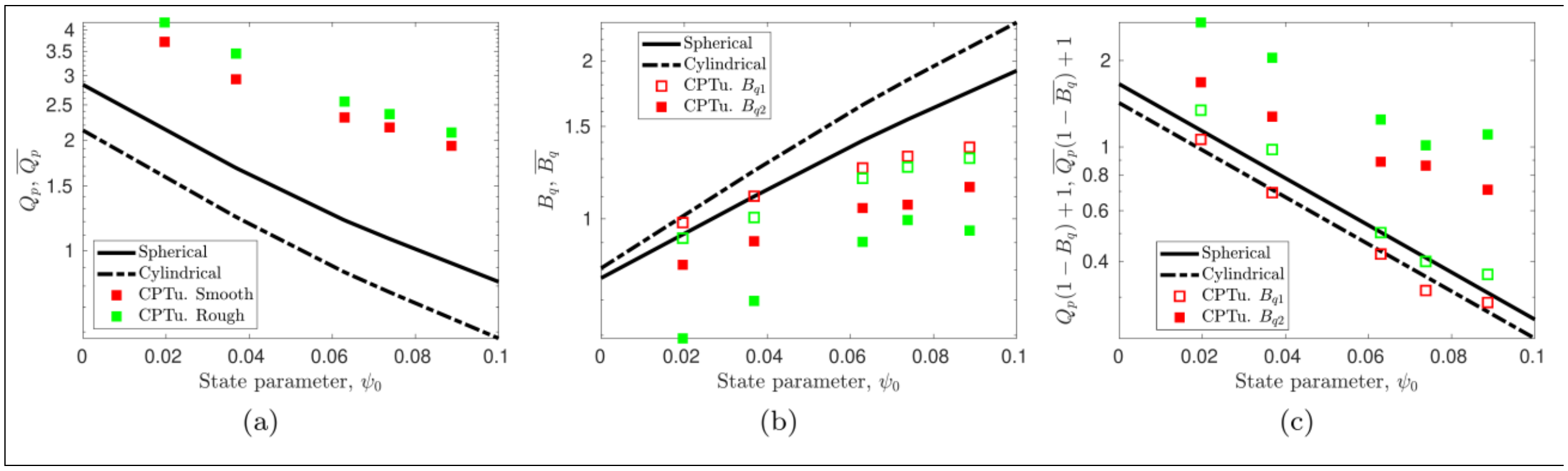

*Figure 9: CPTu simulation results vs cavity expansion reference solutions. Effect of initial state parameter on normalized resistance (a) excess water pressure ratio (b) and normalized effective resistance (c). Comparison of smooth and rough cone.*

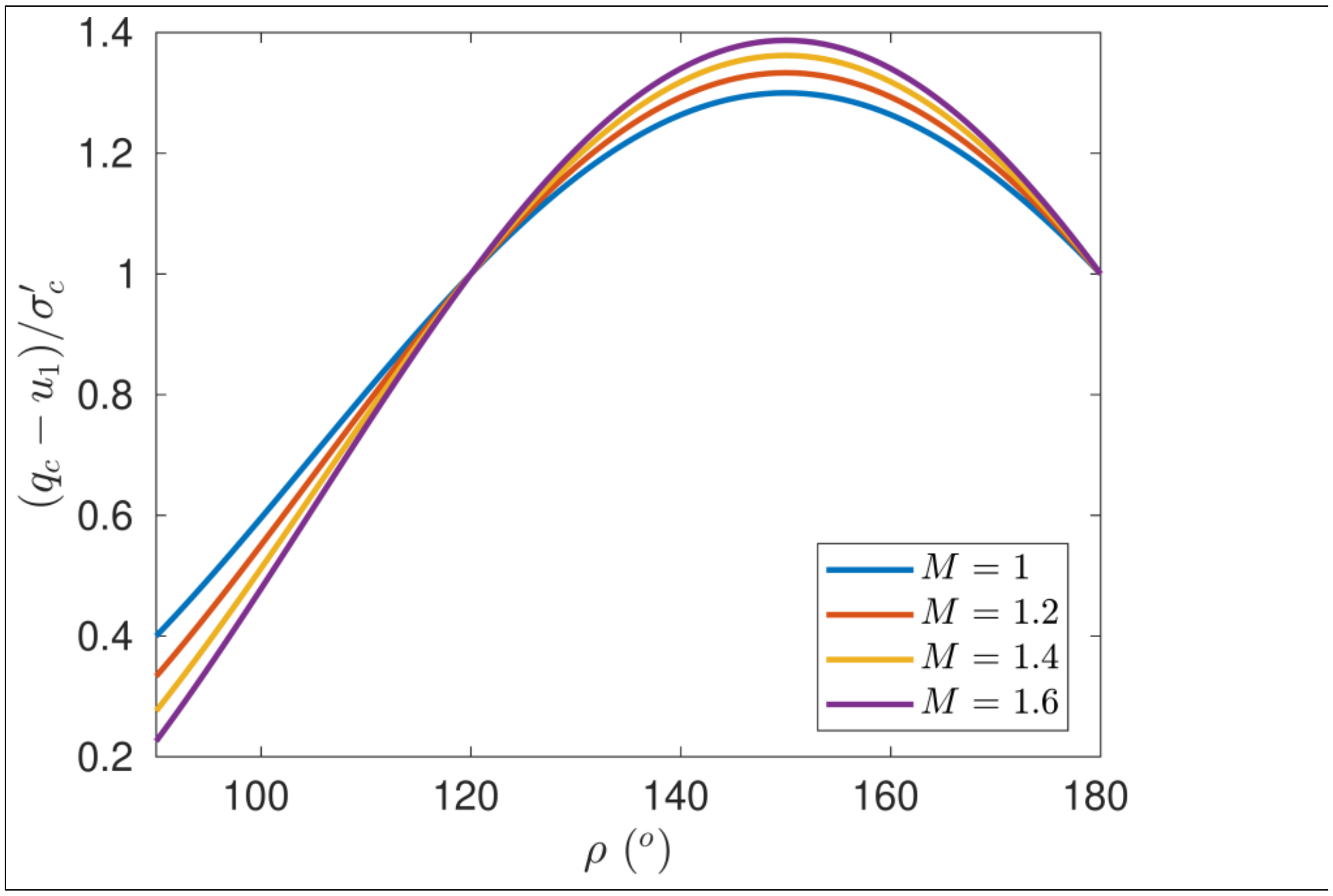

*Figure 10: Ratio between effective tip resistance and cavity resistance in terms of the angle of the major principal effective stress with respect to the vertical(bottom).*

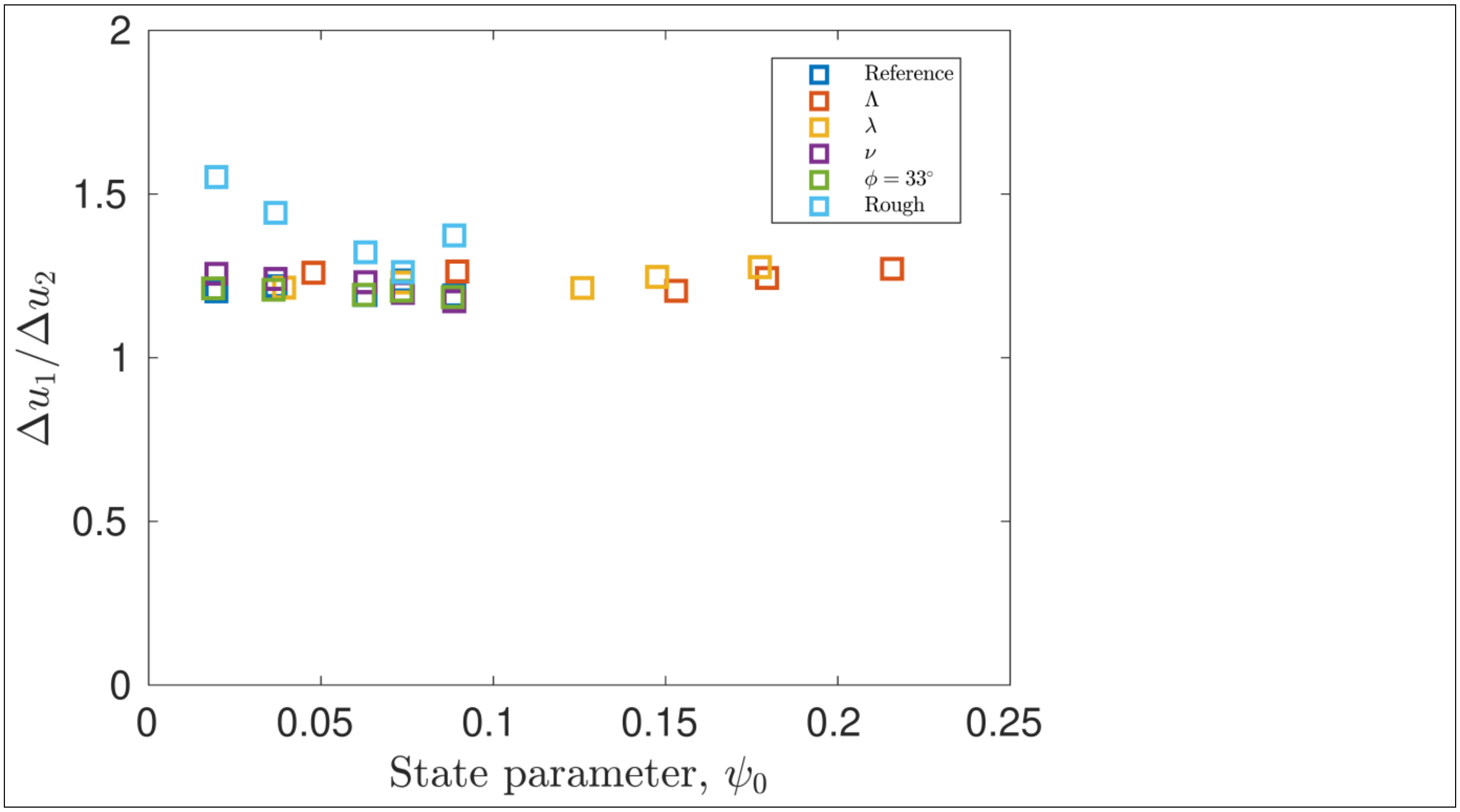

*Figure 11: Ratio β between the u1 and u2 excess water pressure in the G-PFEM simulations*

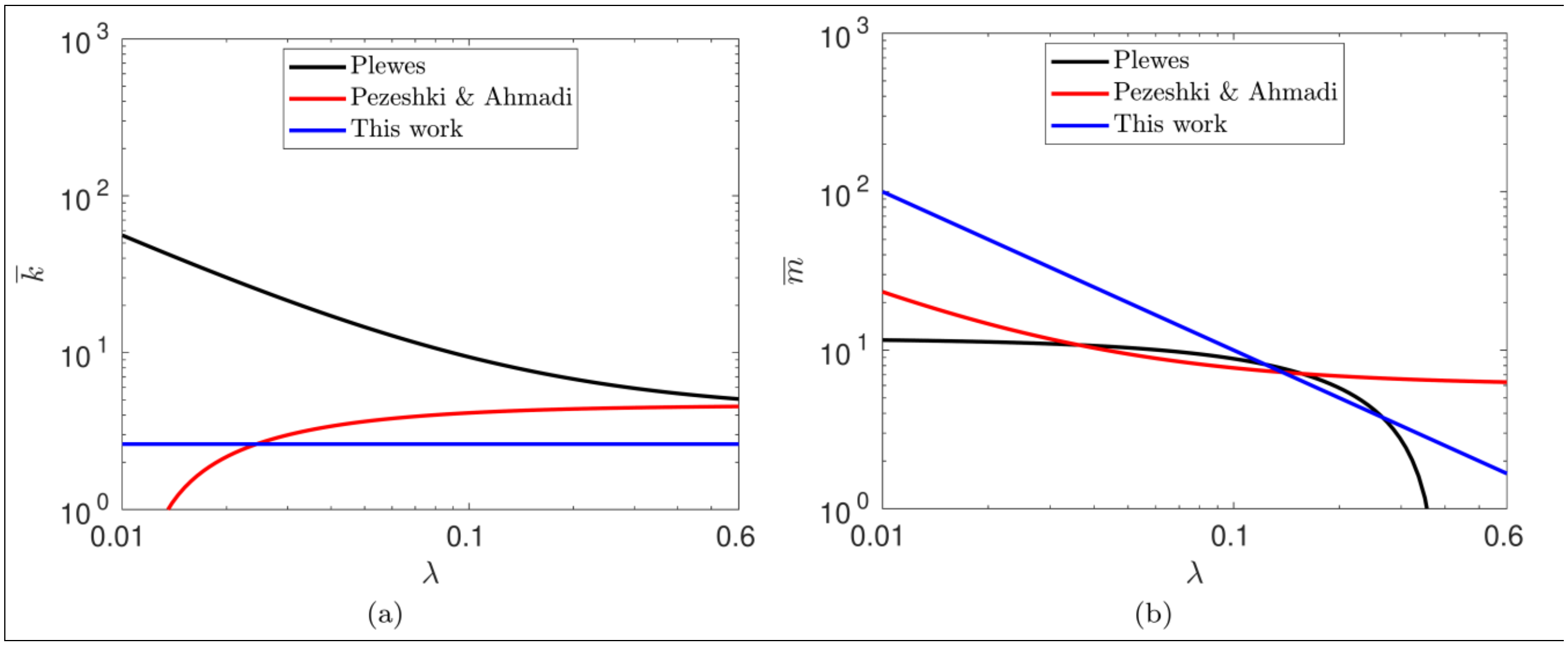

*Figure 12:Effect of $\lambda$ on inversion parameters $\overline{k}$ and $\overline{m}$ for different interpretation methods. $\overline{k}$ plotted for $M = 1.4$. $c_q = 1.35$*

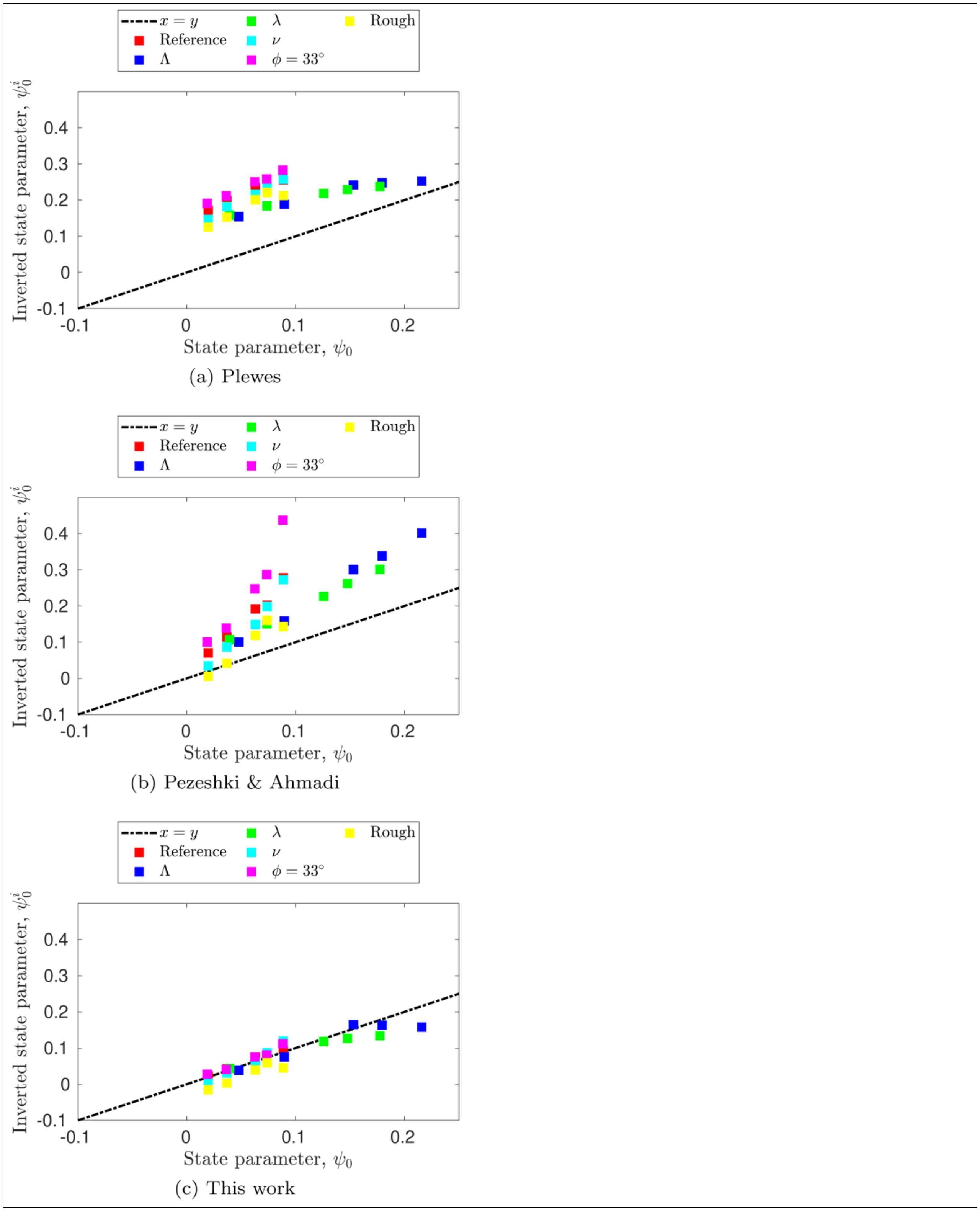

*Figure 13: Comparison of the input state parameter to the numerical simulations and that deduced by different inversion equations. Plewes et al (a), Pezeshki-Ahmadi (b), and this work (c).*

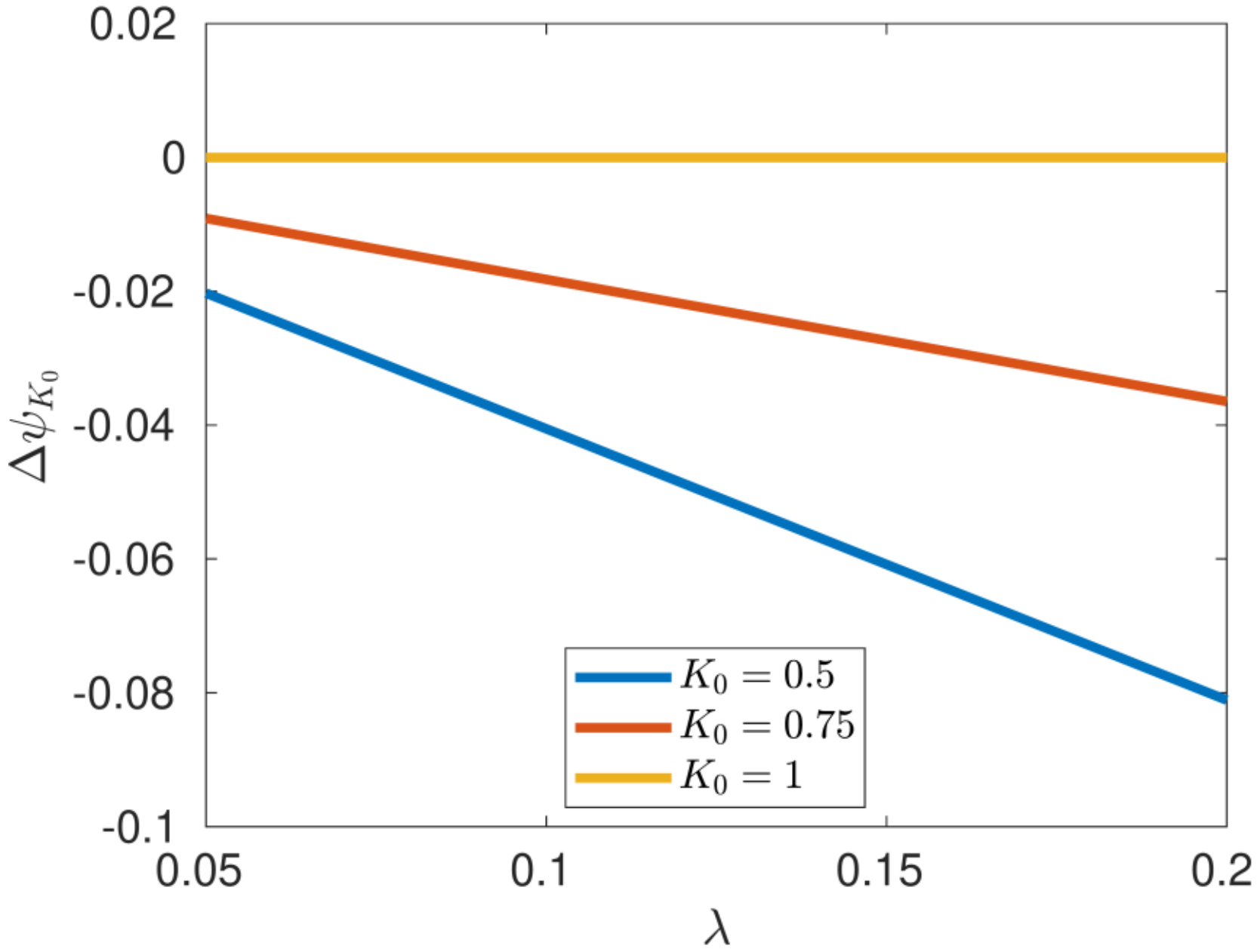

*Figure 14: Effect of the at-rest earth pressure coefficient $K_0$ on the estimated state parameter*

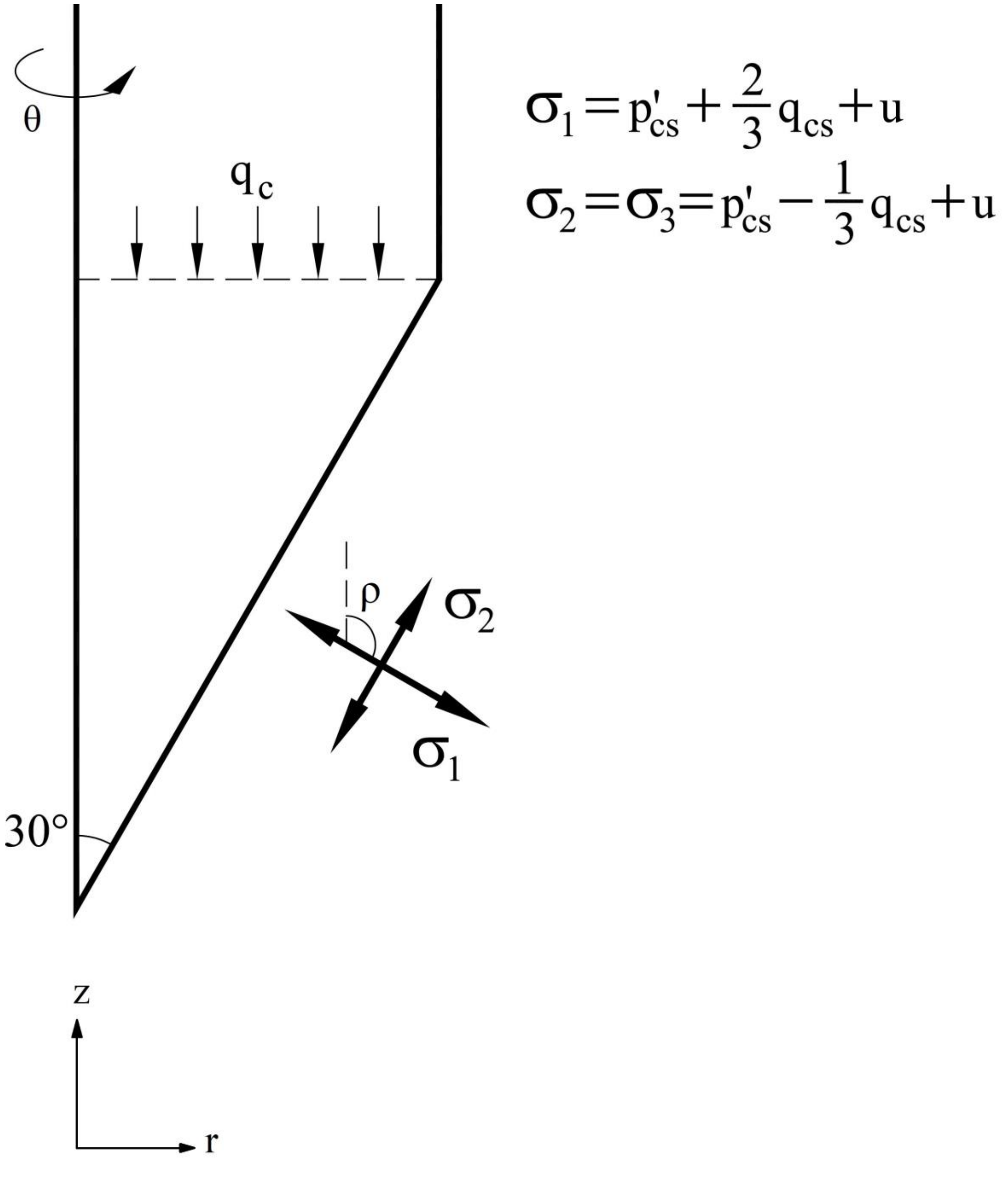

*Figure 15: Assumed stress state near the cone tip for derivation of geometrical term*